\documentclass[DIV=14]{scrartcl}
\usepackage{amssymb}
\usepackage{hyperref}

\usepackage[pdftex]{color}
\usepackage[pdftex]{graphicx}
    \DeclareGraphicsExtensions{.jpg,.pdf,.mps,.png}

\providecommand{\proof}{{\bf Proof: }}
\newtheorem{definition}{DEFINITION}
\newtheorem{lemma}{LEMMA}

\newtheorem{remark}{REMARK}
\newtheorem{example}{EXAMPLE}
\newtheorem{corollary}{COROLLARY}

\newcommand{\eps}{\varepsilon}
\newcommand{\close}{\nolinebreak\mbox{$\blacktriangleleft$}}
\newcommand{\R}{\mathbb{R}}
\newcommand{\N}{\mathbb{N}}
\newcommand{\Z}{\mathbb{Z}}
\newcommand{\spann}{\mathop{\rm span}}

\newenvironment{rmnumerate}
{\begin{enumerate}
}{\end{enumerate}}

\providecommand{\ndofk}{n_{\rm dof}^k}
\providecommand{\ndofK}{n_{\rm dof}^K}

\begin{document}

\title
{Interpolation, projection and hierarchical bases in
discontinuous Galerkin methods}
\author{%
Lutz Angermann\thanks{%
Institut f\"ur Mathematik,
Technische Universit\"at Clausthal,
Erzstra{\ss}e 1, D-38678 Clausthal-Zellerfeld, Germany. 
Email: lutz.angermann@tu-clausthal.de}
\ and
Christian Henke\thanks{%
Institut f\"ur Mathematik,
Technische Universit\"at Clausthal,
Erzstra{\ss}e 1, D-38678 Clausthal-Zellerfeld, Germany.
Email: henke@math.tu-clausthal.de}
}
\date{\today}
\maketitle
{\small
{\bf Abstract:}

The paper presents results on piecewise polynomial approximations
of tensor product type
in Sobolev-Slobodecki spaces
by various interpolation and projection techniques,
on error estimates for quadrature rules and
projection operators based on hierarchical bases, and
on inverse inequalities.
The main focus is directed to applications to discrete conservation laws.

{\bf Keywords:}

Interpolation, quadrature, hierarchical bases, inverse inequalities.

{\bf 2010} {\it Mathematics Subject Classification:}

65\,M\,60, %
65\,D\,32. %
}

\section{Introduction}
The topics covered in this paper belong to the fundamentals
of the analysis of discontinuous Galerkin methods for
partial differential equations.
On the one hand, we have collected and reproved results from areas that
are important for any study of finite element methods,
such as the piecewise polynomial approximation in Sobolev spaces,
quadrature formulas and inverse inequalities (Sections \ref{sec:basics},
\ref{sec:polyappr}, \ref{sec:quadrature}, \ref{sec:inverse}).
On the other hand, we have directed our attention to facts
that are specifically related to particular techniques
such as certain relations between lumping and quadrature effects
or the investigation of fluctuation operators and shock-capturing terms
by means of a hierarchical basis approach
(Sections \ref{sec:quadrature}, \ref{sec:errors}).
A major concern of our study was to trace the dependence of the constants
on both the mesh width and the local polynomial degree.

The paper is organized as follows.
After a brief introduction, which introduces the basic notation,
we investigate polynomial approximations by means of tensor product
elements on affine partitions in Section \ref{sec:polyappr}.
This includes estimates of the reference transformations, which we prove
by the help of a general chain rule. In this way we get a certain insight
into the structure of the occurring constants.
After that we prove error estimates for the Lagrange interpolation
and the $L_2$-projection both with respect to the elements and with respect
to the element edges in the scale of Sobolev-Slobodecki spaces.

In Section \ref{sec:quadrature}, we prepare some important notions such as
quadrature formulas, lumping operators and discrete
$L_2$-projections for later purposes.
In particular, we point out the importance of an suitable choice
of the quadrature points for optimal (w.r.t.\ the local
polynomial degree) error estimates of the Lagrange interpolation.

In the following section we investigate the projection and
interpolation errors for Gauss-Lobatto nodes.
For this purpose we extend the concept of the hierarchical modal basis
to the so-called \emph{embedded hierarchical nodal basis}
and we prove error estimates
for the Lagrange interpolation and for the discrete $L_2$-projection
which are optimal on the elements and almost optimal on the element edges.

As a natural complement to the (direct) estimates from the previous sections,
we present in Section \ref{sec:inverse} inverse inequalities that are based
on generalizations of the Nikolski and Markov inequalities.

\section{Basic notation and definitions}
\label{sec:basics}

Let $\Omega\subset\R^d,$ $d\in\N,$ be a bounded polyhedral domain
with a Lipschitzian boundary (see, e.g., \cite[Def.\ 4.9]{Adams:03}).
$\Omega$ is subdivided by partitions $\mathcal{T}$
(in the sense of \cite[Def.\ 1.49]{Ern:04})
consisting of tensor product elements $T$
(closed as subsets of $\R^d$) with diameter
$h_T := \max_{x,y\in T}\|x - y\|_{\ell^2}.$
Here and in what follows the symbol $\|\cdot\|_{\ell^p},$ $p\in[1,\infty],$
denotes the usual $\ell^p$-norm of (finite) real sequences.
Furthermore the maximal width of a partition $\mathcal{T}$ is defined by
$h := \max\{h_T :\; T\in\mathcal{T}\}.$
To indicate that a particular partition has the maximal width $h,$
we will write $\mathcal{T}_h.$

In this paper, the standard definition of finite elements
$\{T,P_T,\Sigma_T\}$ is used, see e.g.\ \cite[Def.\ 1.23]{Ern:04}.
The finite element space is defined by
$$
W_h := \left\{ w \in L^2(\Omega):\; w|_T \in P_T
\quad
\forall T\in\mathcal{T}_h \right\},
$$
where
$$
P_T\subset W^{l,\infty}(T)\quad\forall T\in\mathcal{T}_h
\quad\mbox{for some}\quad
l\ge 0.
$$
Because of the last requirement,
$W_h$ is a subspace of
$$
W^{l,p}(\mathcal{T}_h)
:= \left\{w \in L^2(\Omega):\; w|_T \in W^{l,p}(T)
\quad\forall T\in\mathcal{T}_h \right\},\quad p\in[1,\infty].
$$
A particular finite element $\{T,P_T,\Sigma_T\}$
is generated by means of a reference element
$\{\hat{T},\hat{P},\hat{\Sigma}\},$
where the geometric reference element $\hat{T}$ is mapped onto
the geometric element $T$ by a $C^1$-diffeomorphism $F_T:\;\hat{T}\to T.$

In the case of a Lagrange finite element $\{\hat{T},\hat{P},\hat{\Sigma}\}$
(in the sense of \cite[Def.\ 1.27]{Ern:04})
with the node set
$\hat\mathcal{N} := \{\hat{x}_1, \ldots , \hat{x}_{\ndofk}\},$
$\ndofk:= \#\hat{\Sigma},$ and the linear forms
$$
\hat{\sigma}_i(\hat{v}) := \hat{v}(\hat{x}_i),\quad 1\le i\le \ndofk,
\quad\forall\hat{v}\in\hat{P} ,
$$
the definitions
\begin{equation}\label{eq:PTdef}
P_T := \{\hat{v}\circ F_T^{-1}:\;\hat{v}\in\hat{P}\}
\end{equation}
and
$$
\sigma_i(v) = \hat{\sigma}_i(v\circ F_T) = \hat{\sigma}_i(\hat{v}),
\quad
v(x)=(\hat{v}\circ F_T^{-1})(x)=\hat{v}(\hat{x}),
\quad
x = F_T(\hat{x}),
\quad
1\le i\le \ndofk,
$$
are used.
The nodal basis of $P_T$ is obtained by an analogous transformation
of the reference shape functions
$$
\{\hat{\varphi}_1, \ldots ,\hat{\varphi}_{\ndofk}\}
\quad\mbox{with}\quad
\hat{\varphi}_i(\hat{x}_j) = \delta_{ij},
\quad 1\le i,j\le \ndofk.
$$
\begin{definition}
A partition $\mathcal{T}_h$ is called \emph{affine} if the mapping
$F_T$ is affine for all $T\in\mathcal{T}_h,$ i.e. if
$F_T(\hat{x})= J_T \hat{x} + b_T$ with $J_T \in\R^{d,d},$ $\det J_T \ne 0,$
$b_T \in\R^d.$
\end{definition}
In addition, the following properties of $\mathcal{T}_h$ are important.
\begin{definition}%
[locally quasiuniform, shape regular]
A family of affine partitions $\{\mathcal{T}_h\}_{h>0}$ is called
\emph{locally quasiuniform} if there exists a constant $\sigma_0>0$
such that
\begin{equation}\label{lokalquasiuniform} %
\forall h>0\ \forall T\in\mathcal{T}_h:\ \sigma_T := \frac{h_T}{\varrho_T}
\le\sigma_0,
\end{equation}
where $\varrho_T$ denotes the diameter of the largest ball contained in $T.$
\end{definition}
\begin{definition}%
[quasiuniform]
A family of partitions $\{\mathcal{T}_h\}_{h>0}$ is called
\emph{quasiuniform} if it is locally quasiuniform and if a constant
$C_{qu} > 0$ exists with
$$
\forall h>0\ \forall T\in\mathcal{T}_h:\ h_T \ge C_{qu}h.
$$
\end{definition}

\section{Polynomial approximation using tensor product elements}
\label{sec:polyappr}
In what follows we will investigate thoroughly the approximation
of functions on affine partitions.
We set
$\hat{T} := I^d,$ $I := [-1,1]$ and
$\hat{P} := \mathbb{Q}_k(\hat{T})$ with
$$
\mathbb{Q}_k(\hat{T})
:=\spann_{\alpha\in\N_0^d,\,\|\alpha\|_{\ell^\infty}\le k}\{\hat{x}^\alpha\},
\quad\hat{x}\in\hat{T},\ k \in\N_0:=\N\cup\{0\}.
$$
A function
$\hat{\varphi}_\alpha\in\{\hat{\varphi}_1, \ldots ,\hat{\varphi}_{\ndofk}\},$
$\alpha\in\N_0^d,$
can be written as a product of univariate Lagrange polynomials.
Namely, denote by
$\{\hat{\varphi}_0^k,\hat{\varphi}_1^k, \ldots ,\hat{\varphi}_k^k\}$
a basis of the space of univariate polynomials of maximum degree $k.$
Then, for any multiindex $\alpha\in\N_0^d$ with $\|\alpha\|_{\ell^\infty}\le k,$
$$
\hat{\varphi}_\alpha(\hat{x})=\hat{\varphi}_{\alpha_1}^k(\hat{x}_1)
\hat{\varphi}_{\alpha_2}^k(\hat{x}_2)
\cdots\hat{\varphi}_{\alpha_d}^k(\hat{x}_d).
$$
The space $\mathbb{Q}_k(T)$ is defined according to (\ref{eq:PTdef}).
Similarly, based on the definition
$$
\mathbb{P}_k(\hat{T})
:=\spann_{\alpha\in\N_0^d,\,\|\alpha\|_{\ell^1}\le k}\{\hat{x}^\alpha\},
\quad\hat{x}\in\hat{T},\ k \in\N_0,
$$
the space $\mathbb{P}_k(T)$ of polynomials of maximum degree $k$ can
be introduced.

As a consequence of the affine structure of the transformation
between the reference element $\hat{T}$ and the element $T$
we can prove the following estimates.
\begin{lemma} %
\label{lemma_affin_trans}
For $l \ge 0$ and $1 \le p \le \infty,$ $1/\infty:=0,$
there exists a constant $C_{l,d}\ge 1$
such that, for $T \in \mathcal{T}_h,$ $\mathcal{T}_h$ affine, and
$w \in W^{l,p}(T),$ $\hat{w}=w\circ F_T,$ the following estimates hold:
\begin{eqnarray}
| \hat{w}|_{l,p,\hat{T}} &\le & C_{l,d} \|J_T\|_{\ell^2}^l |\det J_T|^{-1/p}
|w|_{l,p,T} ,
\label{trans_ungl1} \\
| w|_{l,p,T} &\le & C_{l,d} \|J_T^{-1}\|_{\ell^2}^l |\det J_T|^{1/p}
|\hat{w}|_{l,p,\hat{T}} ,
\label{trans_ungl2}
\end{eqnarray}
where $C_{l,d}$ depends only on $l$ and $d.$ In particular, $C_{0,d}=1.$
\end{lemma}
\proof
By Fa\`{a} di Bruno's formula (see, e.g. \cite{Johnson:02}),
for $\hat{w}=w \circ F_T, \, |\alpha|=l,$ $P:=\N_0^d \setminus \{0\}$ and
$$
\mathcal{A}
:=\left\{ a:P \to \N_0^d : \sum_{\gamma \in P} a(\gamma) =\beta
\mbox{ and } \sum_{\gamma \in P} |a(\gamma)| \gamma =\alpha \right\}
$$
we have the representation
$$
\begin{array}{rcl}
\hat{\partial}^\alpha (w \circ F_T)(\hat{x})
&=&\alpha! \sum_{|\beta| \le |\alpha|} (\partial^\beta w)(F_T(\hat{x}))
\sum_{a \in \mathcal{A}}
\prod_{\gamma \in P} \frac{1}{a(\gamma)!} \left[ \frac{(\hat{\partial}^\gamma
F_T)(\hat{x})}{\gamma!} \right]^{a(\gamma)} \\
&=&\alpha! \sum_{|\beta| \le |\alpha|} (\partial^\beta w)(F_T(\hat{x}))
\sum_{a \in \mathcal{A}}
\prod_{j = 1}^{d} \frac{1}{a(e_j)!} \left[ (\hat{\partial}_j F_T)
\right]^{a(e_j)} \\
&=&\alpha! \sum_{|\beta| =l} (\partial^\beta w)(F_T(\hat{x}))
\sum_{a \in \mathcal{A}}
\prod_{j = 1}^{d} \frac{1}{a(e_j)!} \left[ (\hat{\partial}_j F_T)
\right]^{a(e_j)},
\end{array}
$$
because the conditions
$l=|\alpha|=|\sum_{\gamma \in \N^m} |a(\gamma) |
\gamma|=|\sum_{j=1}^d |a(e_j)| e_j|= \sum_{j=1}^d |a(e_j)|$ and
$|\beta|=|\sum_{j=1}^d a(e_j)|=\sum_{j=1}^d |a(e_j)|$
imply that
$|\alpha|=|\beta|=l$ if ``empty'' sums are neglected.
The absolute value of the left-hand side can be estimated as
$$
\begin{array}{rcl}
|\hat{\partial}^\alpha (w \circ F_T)(\hat{x})|
&\le& \alpha! \sum_{|\beta| =l} \left| (\partial^\beta w)(F_T(\hat{x})) \right|
\sum_{a \in \mathcal{A}}
\prod_{j = 1}^{d} \frac{1}{a(e_j)!} \left| (\hat{\partial}_j F_T)^{a(e_j)}
\right| \\
&\le& \alpha! \sum_{|\beta| =l} \left| (\partial^\beta w)(F_T(\hat{x})) \right|
\sum_{a \in \mathcal{A}}
\prod_{j = 1}^{d} \frac{1}{a(e_j)!} \left\| \hat{\partial}_j
F_T\right\|_{\ell^\infty} ^{|a(e_j)|} \\
&\le& \alpha! \sum_{|\beta| =l} \left| (\partial^\beta w)(F_T(\hat{x})) \right|
\sum_{a \in \mathcal{A}}
\prod_{j = 1}^{d} \frac{1}{a(e_j)!} \left\|J_T \right\|_{\max} ^{l} \\
&=& C_{\alpha,d} \sum_{|\beta| =l} \left| (\partial^\beta w)(F_T(\hat{x}))
\right|
\left\|J_T \right\|_{\max} ^{l} \\
&\le& C_{\alpha,d} \sum_{|\beta| =l} \left| (\partial^\beta w)(F_T(\hat{x}))
\right|
\left\|J_T \right\|_{\ell^2} ^{l}.
\end{array}
$$
For $p<\infty,$
using H\"older's inequality for sums and observing that
$\sum_{|\beta|=l}1={d+l-1 \choose l},$
on the reference element we have that
$$
\| \hat{\partial}^{\alpha} \hat{w}\|_{0,p,\hat{T}}^p
\le {d+l-1 \choose l}^{(p-1)} C_{\alpha,d}^p
\|J_T\|_{\ell^2}^{lp} \sum_{|\beta|=l} \|\partial^\beta w \circ F_T
\|_{0,p,\hat{T}}^p,
$$
Applying the substitution rule, we get
$$
\| \hat{\partial}^\alpha \hat{w}\|_{0,p,\hat{T}}^p \le {d+l-1 \choose l}^{(p-1)}
C_{\alpha,d}^p \|J_T\|_{\ell^2}^{lp} |\det J_T|^{-1} |w|_{l,p,T}^p.
$$
Summing up w.r.t.\ $\alpha,$ the estimate
$$
\begin{array}{rcl}
|\hat{w}|_{l,p,\hat{T}}^p&=&\sum_{|\alpha|=l} \|\hat{\partial}^\alpha
\hat{w}\|_{0,p,\hat{T}}^p \\
&\le& \sum_{|\alpha|=l} {d+l-1 \choose l}^{(p-1)} C_{\alpha,d}^p
\|J_T\|_{\ell^2}^{lp} |\det J_T|^{-1}|w|_{l,p,T}^p \\
&\le& \left(\max\limits_{|\alpha|=l}\left( C_{\alpha,d}\right) \right)^p
{d+l-1 \choose l}^{p} \|J_T\|_{\ell^2}^{lp} |\det J_T|^{-1}|w|_{l,p,T}^p
\end{array}
$$
proves the statement with
$C_{l,d}:={d+l-1 \choose l} \max\limits_{|\alpha|=l}\left( C_{\alpha,d} \right).$
For $p=\infty$ we see that
$$
\begin{array}{rcl}
|\hat{w}|_{l,\infty,\hat{T}} &=& \max\limits_{|\alpha|=l}
\|\hat{\partial}^\alpha \hat{w}\|_{0,\infty,\hat{T}} \\
&\le& \max\limits_{|\alpha|=l} \left( C_{\alpha,d}\right) \|J_T\|_{\ell^2}^l
\|\sum_{|\beta|=l} \left| \left( \partial^\beta w \right) \left( F_T(\hat{x})
\right) \right| \|_{0,\infty,\hat{T}} \\
&=& \max\limits_{|\alpha|=l} \left( C_{\alpha,d}\right) \|J_T\|_{\ell^2}^l
\|\sum_{|\beta|=l} \left| \left( \partial^\beta w \right) \left( x \right)
\right| \|_{0,\infty,T}\\
&\le& \max\limits_{|\alpha|=l} \left(C_{\alpha,d}\right) \|J_T\|_{\ell^2}^l
{d+l-1 \choose l} \| \max\limits_{|\beta|=l} \left| \partial^\beta w
\right| \|_{0,\infty,T}\\
&=& C_{l,d} \|J_T\|_{\ell^2}^l \max\limits_{|\beta|=l} \| \partial^\beta w
\|_{0,\infty,T}
\end{array}.
$$
Since $F_T:\hat{T} \to T$ is bijective, the second estimate follows obviously.
\close
\begin{lemma} %
The following estimates are valid:
$$
|\det J_T|=\frac{|T|}{|\hat{T}|},\quad \|J_T\|_{\ell^2} \le
\frac{h_T}{\rho_{\hat{T}}}
\quad\mbox{and}\quad
\|J_T^{-1}\|_{\ell^2} \le
\frac{h_{\hat{T}}}{\rho_T}.
\label{help_scaling}
$$
\end{lemma}
\proof
The first relation is classical. The proof of the inequalities
is easy, see, e.g., \cite[Thm.\ 15.2]{Ciarlet:91a}.
\close

\begin{corollary} %
Given a locally quasiuniform family $\{\mathcal{T}_h \}_{h >0}$
of affine partitions with the reference element $\hat{T}=I^d.$
Then:
\begin{equation}\label{help_scaling2}
\|J_T\|_{\ell^2}
\le \frac{h_T}{2}, \quad \|J_T^{-1}\|_{\ell^2}
\le \frac{2\sqrt{d}\sigma_0}{h_T}
\end{equation}
and
\begin{equation}\label{help_scaling3}
\begin{array}{c}\displaystyle
\det J_T=\prod_{i=1}^d \lambda_i^{1/2}\le\lambda_{\max}^{d/2}
=\|J_T\|_{\ell^2}^d\le 2^{-d} h_T^{d},
\qquad %
|T| =|\hat{T}| |\det J_T| \le h_T^d,
\end{array}
\end{equation}
where $\lambda_i,$ $1 \le i \le d,$
are the eigenvalues of the matrix $J_T^\top J_T.$
\end{corollary}
Now we are ready to define the Lagrange interpolation operator as follows:
$$
I_h^k:\; C^0(T) \ni v \mapsto I_h^k v:=\sum_{i=1}^{\ndofk} v(x_i)
\varphi_i \in \mathbb{Q}_k(T) \subset W^{l,\infty}(T), \quad x_i \in
\mathcal{N}_T:=F_T(\hat\mathcal{N}).
$$
We mention that the results from this section on the Lagrange
interpolation operator do not impose any conditions
w.r.t.\ the location of the nodes.
Later we will formulate statements about the interpolation operator
which use the special properties of Gauss-Lobatto quadrature points.

Simple calculations show that
\begin{eqnarray}
I_h^k v&=&v \quad\forall v \in \mathbb{Q}_k(T),
\label{int_eig_inv} \\
\|I_h^k v\|_{l,p,T} &\le & C \|v\|_{0,\infty,T},
\quad l \ge 0,\ p \in[1,\infty],
\label{int_eig_stetig}
\end{eqnarray}
where $C=C(T,\{\varphi_i\}_{i=1}^{\ndofk},\ndofk,l,p)>0.$
Furthermore we have the following error estimates.
\begin{lemma}\label{interpolationsfehler} %
Let $T$ be an element of an affine partition $\mathcal{T}_h$
such that the corresponding family of partitions $\{\mathcal{T}_h\}_{h>0}$
is locally quasiuniform.
Assume that
$1\le l \le k+1,$ $l \in \N,$ $p \in [1,\infty]$ and $lp > d.$
Then, for the Lagrange interpolation operator $I_h^k,$
there exist constants $C > 0$ independent of $h_T$ such that
\begin{eqnarray*}
|v - I_h^k v|_{r,p,T} & \le & C h_T^{l-r} |v|_{l,p,T}, \quad 0\le r \le l,\\
|v - I_h^k v|_{r,p,E} & \le & C h_T^{l-1/p-r} |v|_{l,p,T}, \quad 0\le r \le l-1/p,
\end{eqnarray*}
for all $v \in W^{l,p}(T),$
where $E$ denotes a face of $T.$
\end{lemma}
The proof relies on the following interpolation inequality.
\begin{lemma}[interpolation inequality]\label{l:interpolineq2}
Let $l_0,l_1 \in \N_0,\, l_0 \ne l_1,$ $1 \le p \le \infty$
and $G$ be a bounded domain with a Lipschitzian boundary.
Define, for $0 < \theta < 1,$
$l_{\theta}:=(1-\theta)l_0 + \theta l_1.$
Then
\begin{equation}\label{interpolationsungleichung2}
\forall u \in W^{l_0,p}(G) \cap W^{l_1,p}(G):
\quad
\|u\|_{l_{\theta},p,G} \le C \|u\|_{l_0,p,G}^{1-\theta}\|u\|_{l_1,p,G}^\theta \,.
\end{equation}
\end{lemma}
\textbf{Proof} of Lemma \ref{l:interpolineq2}:
We first mention that the statement for the case $G=\R^d$
is a consequence of \cite[1.3.3 (g)]{Triebel:78} and
\cite[Def.\ 6.2.2, Thm.\ 6.2.3, Thm.\ 6.2.4 and Thm.\ 6.4.5 (3),(4)]{Bergh:76}.
If $G$ is a bounded domain with a Lipschitzian boundary,
then there exists a total extension operator $E_G$
(see \cite[Thm.\ 5.24]{Adams:03} or \cite[Ch.\ 6, Thm.\ 5]{Stein:70})
such that
$$
\begin{array}{rcl@{\quad}l}
\|u\|_{l_{\theta},p,G}
& = & \|E_Gu\|_{l_{\theta},p,G}
&\mbox{(by \cite[(i) in (5.17)]{Adams:03})}\\
&\le& \|E_Gu\|_{l_{\theta},p,\R^d} \\
&\le& C \|E_Gu\|_{l_0,p,\R^d}^{1-\theta}
\|E_Gu\|_{l_1,p,\R^d}^{\theta}\\
&\le& C \|u\|_{l_0,p,G}^{1-\theta}\|u\|_{l_1,p,G}^{\theta}
&\mbox{(by \cite[(ii) in (5.17)]{Adams:03})}.
\qquad\close
\end{array}
$$
\begin{remark}
Using further results on Stein's extension operator
\cite[p.\ 185 and Thm.\ 1]{Kalyabin:85}, the interpolation inequality
(\ref{interpolationsungleichung2}) can be extended to the
parameter set $l_0,l_1 \in \R_+,\, l_0 \ne l_1,$ $1 \le p_0,p_1 \le \infty,$
as follows, where
$\frac{1}{p_{\theta}} := \frac{1-\theta}{p_0} + \frac{\theta}{p_1}$:
$$
\forall u \in W^{l_0,p_0}(G) \cap W^{l_1,p_1}(G):
\quad
\|u\|_{l_{\theta},p_{\theta},G}
\le C \|u\|_{l_0,p_0,G}^{1-\theta}\|u\|_{l_1,p_1,G}^\theta \,.
$$
\end{remark}
\textbf{Proof} of Lemma \ref{interpolationsfehler}:
By the triangle inequality, (\ref{int_eig_stetig}), and the embedding
theorem \cite[Thm.\ 4.12, Part II]{Adams:03},
for $1\le l\le k+1,$ $r \in \N_0,$ $r\le l,$ $lp>d$ we have that
\begin{equation}\label{error_element_norm}
\begin{array}{rcl}
\|\hat{v} - \hat{I}_h^k \hat{v}\|_{r,p,\hat{T}}
&\le& \|\hat{v}\|_{r,p,\hat{T}} + \|\hat{I}_h^k \hat{v}\|_{r,p,\hat{T}}
\le \|\hat{v}\|_{l,p,\hat{T}} + C\|\hat{v}\|_{0,\infty,\hat{T}} \\
&\le& C \|\hat{v}\|_{l,p,\hat{T}} \quad \forall \hat{v} \in W^{l,p}(\hat{T}).
\end{array}
\end{equation}
Furthermore, by (\ref{int_eig_inv}) and
$\mathbb{P}_k(\hat{T}) \subset \mathbb{Q}_k(\hat{T})$,
\begin{eqnarray*}
|\hat{v}-\hat{I}_h^k \hat{v}|_{r,p,\hat{T}}
&\le& \|\hat{v}-\hat{I}_h^k \hat{v}\|_{r,p,\hat{T}}
=\inf\limits_{\hat{p} \in \mathbb{P}_{k}(\hat{T})}
\|(Id-\hat{I}_h^k)(\hat{v}+\hat{p})\|_{r,p,\hat{T}} \\
&\le& \inf\limits_{\hat{p} \in \mathbb{P}_{l-1}(\hat{T})}
\|(Id-\hat{I}_h^k)(\hat{v}+\hat{p})\|_{r,p,\hat{T}},
\end{eqnarray*}
where $Id$ denotes the identity operator.
From (\ref{error_element_norm}) we see that
$$
|\hat{v}-\hat{I}_h^k \hat{v}|_{r,p,\hat{T}}
\le C \inf\limits_{\hat{p} \in \mathbb{P}_{l-1}(\hat{T})}
\|\hat{v}+\hat{p}\|_{l,p,\hat{T}} \le C |\hat{v}|_{l,p,\hat{T}},
$$
where the last estimate is a consequence of the Deny-Lions lemma
(\cite[Thm.\ 14.1]{Ciarlet:91a}).
Note that the constant $C$ depends on the parameters of the reference
element.
The application of Lemma \ref{lemma_affin_trans} results in an estimate
on the element $T$:
$$
\begin{array}{rcl}
|v-I_h^k v|_{r,p,T}
&\le& C\|J_T^{-1}\|_{\ell^2}^r |\det J_T|^{1/p}|\hat{v}-\hat{I}_h^k\hat{v}|_{r,p,\hat{T}}\\
&\le& C\|J_T^{-1}\|_{\ell^2}^r |\det J_T|^{1/p}|\hat{v}|_{l,p,\hat{T}}\\
&\le& C\left(\|J_T\|_{\ell^2} \|J_T^{-1}\|_{\ell^2} \right)^r \|J_T\|_{\ell^2}^{l-r}
|v|_{l,p,T}\\
&\le& C \left( \frac{h_T}{\rho_T}\right)^r h_T^{l-r} |v|_{l,p,T}
\le C h_T^{l-r} |v|_{l,p,T}.
\end{array}
$$
where we have used the condition (\ref{lokalquasiuniform}) in the last step.

In the case $r\in\R_+\setminus\N_0,$ we apply the interpolation inequality
(\ref{interpolationsungleichung2}) with
$l_0:=\lfloor r \rfloor:=\max\limits_{z \in \Z,\, \{z\} \le r} z,$
$l_1:=\lceil r \rceil:=\min\limits_{z \in \Z,\, z \ge r} \{z\}$
and $\theta: = r-\lfloor r\rfloor$:
\begin{eqnarray*}
|v-I_h^k v|_{r,p,T} &\le & C |v-I_h^k v|_{\lfloor r \rfloor, p,T}^{1-\theta}
|v-I_h^k v|_{\lceil r \rceil,p,T}^\theta \\
&\le & C h_T^{l-r} |v|_{l,p,T} ,\quad \lceil r \rceil \le l.
\end{eqnarray*}
The proof of the second estimate %
runs similarly.

For $1 \le l \le k+1,$ $p\in (1,\infty),$ and any face $\hat{E}$
of the reference element $\hat{T}$ with $E=F_T \hat{E}$
we have, by the special trace theorem for the faces of $\hat{T}$
(see \cite[Thm.\ 2.5.4]{Necas:67}) and (\ref{error_element_norm}), that
\begin{equation}\label{error_kante_norm1}
\|\hat{v}-\hat{I}_h^k \hat{v}\|_{l-1/p,p,\hat{E}}
\le C\|\hat{v}- \hat{I}_h^k \hat{v}\|_{l,p,\hat{T}}
\le C \|\hat{v}\|_{l,p,\hat{T}}
\quad\forall \hat{v} \in W^{l,p}(\hat{T}).
\end{equation}
Furthermore, \cite[Lemma 2.5.4]{Necas:67} implies that, for any
$r\in[0,l-1/p),$
$$
\|\hat{v}-\hat{I}_h^k \hat{v}\|_{r,p,\hat{E}}
\le C \|\hat{v}-\hat{I}_h^k \hat{v}\|_{l-1/p,p,\hat{E}}
\quad\forall \hat{v} \in W^{l,p}(\hat{T}).
$$
Combining this estimate with (\ref{error_kante_norm1}) we arrive,
for $p\in (1,\infty)$ and $r\in[0,l-1/p],$ at
\begin{equation}\label{error_kante_norm2}
\|\hat{v}-\hat{I}_h^k \hat{v}\|_{r,p,\hat{E}}
\le C \|\hat{v}\|_{l,p,\hat{T}}
\quad\forall \hat{v} \in W^{l,p}(\hat{T}).
\end{equation}
In the case $p=1,$ we conclude from \cite[Thm.\ 1.II]{Gagliardo:57}
by a similar argument as in \cite[Lemma 2.5.4]{Necas:67} that
$$
\|\hat{v}-\hat{I}_h^k \hat{v}\|_{l-1,1,\hat{E}}
\le C\|\hat{v}- \hat{I}_h^k \hat{v}\|_{l,1,\hat{T}}
\quad\forall \hat{v} \in W^{l,1}(\hat{T}),
$$
hence,
$$
\|\hat{v}-\hat{I}_h^k \hat{v}\|_{r,1,\hat{E}}
\le C \|\hat{v}\|_{l,1,\hat{T}}
\quad\forall \hat{v} \in W^{l,1}(\hat{T})
$$
for all $r\in[0,l-1].$
The case $p=\infty$ is a simple consequence of the fact that
the trace operator is the classical restriction
due to the embedding theorem \cite[Thm.\ 2.3.8]{Necas:67}:
$$
\|\hat{v}-\hat{I}_h^k \hat{v}\|_{r,\infty,\hat{E}}
\le C \|\hat{v}\|_{l,\infty,\hat{T}}
\quad\forall \hat{v} \in W^{l,\infty}(\hat{T})
$$
for all $r\in[0,l].$
Thus (\ref{error_kante_norm2})
is proved for all $p \in [1,\infty]$ with $lp > d$ and all $r\in[0,l-1/p].$

Now, if additionally $r\in\N_0,$ we have that
$$
\begin{array}{rcl}
|\hat{v}-\hat{I}_h^k \hat{v}|_{r,p,\hat{E}}
&\le& \|\hat{v}-\hat{I}_h^k \hat{v}\|_{r,p,\hat{E}}
=\inf\limits_{\hat{p} \in \mathbb{P}_{k}(\hat{T})} \|(Id-\hat{I}_h^k)(\hat{v}
+\hat{p})\|_{r,p,\hat{E}} \\
&\le& \inf\limits_{\hat{p} \in \mathbb{P}_{l-1}(\hat{T})}
\|(Id-\hat{I}_h^k)(\hat{v}+\hat{p})\|_{r,p,\hat{E}} \\
&\le& C\inf\limits_{\hat{p} \in \mathbb{P}_{l-1}(\hat{T})}
\|\hat{v} +\hat{p}\|_{l,p,\hat{T}} \le C |\hat{v}|_{l,p,\hat{T}}.
\end{array}
$$
Using the interpolation inequality (\ref{interpolationsungleichung2})
and performing the back-transformation, we get
$$
\begin{array}{rcl}
|v-I_h^k v|_{r,p,E}
&\le& C |v-I_h^k v|_{\lfloor r \rfloor,p,E}^{1-\theta}
|v -I_h^k v|_{\lceil r \rceil,p,E}^{\theta}\\
&\le& C \|J_E^{-1}\|_{\ell^2}^r |\det J_E|^{1/p} |\hat{v}|_{l,p,\hat{T}}\\
&\le& C \left(\|J_T\|_{\ell^2} \|J_E^{-1}\|_{\ell^2} \right)^r \|J_T\|_{\ell^2}^{l-r}
\left( \frac{|E| |\hat{T}|}{|\hat{E}||T|} \right)^{1/p} |v|_{l,p,T}\\
&\le& C \left( \frac{h_T}{\rho_E}\right)^r h_T^{l-r} h_T^{-1/p} |v|_{l,p,T}\\
&\le& C h_T^{l-1/p-r} |v|_{l,p,T},
\end{array}
$$
where we have used the simple estimate $\rho_T\le \rho_E$
together with (\ref{lokalquasiuniform}).
\close

\begin{corollary}\label{Folg_interpolationsfehler} %
Let $T$ be an element of an affine partition $\mathcal{T}_h$
such that the corresponding family of partitions $\{\mathcal{T}_h\}_{h>0}$
is locally quasiuniform.
Assume that $l \in\N,$ $p \in [1,\infty]$ and $lp > d.$
Then, for the Lagrange interpolation operator $I_h^k,$
there exist constants $C > 0$ independent of $h_T$ such that
\begin{eqnarray*}
|v - I_h^k v|_{r,p,T}
&\le& C h_T^{\min\{k+1,l\}-r} \|v\|_{l,p,T}, \quad 0\le r \le l,\\
|v - I_h^k v|_{r,p,E}
&\le& C h_T^{\min\{k+1,l\}-1/p-r} \|v\|_{l,p,T}, \quad 0\le r \le l-1/p,
\end{eqnarray*}
for all $v \in W^{l,p}(T).$
\end{corollary}
\proof
For $1\le l \le k+1,$ the statement coincides with Lemma \ref{interpolationsfehler}.
In the case $k+1<l,$ $r\in\R_+,$ Lemma \ref{interpolationsfehler} implies that
$$
|v-I_h^k v|_{r,p,T} \le C h_T^{k+1-r} |v|_{k+1,p,T}.
$$
From $h_T^{k+1-r}|v|_{k+1,p,T} \le C h_T^{\min\{k+1,l\}-r} \|v\|_{l,p,T}$
the first estimate follows. %
The proof of the second estimate %
runs analogously.
\close

A further possibility of approximating functions in Sobolev spaces
is given by the projection w.r.t.\ the $L^2$ inner product.
\begin{definition} %
The orthogonal $L^2$-projection $P_h^k:\;L^2(T) \to\mathbb{Q}_k(T)$
is defined, for $v\in L^2(T),$ by
\begin{equation}\label{L2_proj}
(v -P_h^k v, w)_{0,T}=0\quad \forall w \in \mathbb{Q}_k(T).
\end{equation}
\end{definition}
In general, the relation (\ref{L2_proj}) is equivalent to a system of
linear algebraic equations. It can be solved easily provided an
$L^2$-orthogonal basis is used. For elements $T$ of an affine partition
and any multiindex $\alpha\in\N_0^d,$ $\|\alpha\|_{\ell^\infty}\le k,$
we have that
$$
\int_T \psi_\alpha \psi_\beta \, dx = |\det J_T|\int_{\hat{T}}
\hat{\psi}_\alpha \hat{\psi}_{\beta} \, d\hat{x}=|\det J_T| \prod_{i=0}^d
\frac{1}{2 \alpha_i +1} \delta_{\alpha \beta}=\rho_\alpha^J
\delta_{\alpha,\beta},
$$
where
$$
\hat{\psi}_{\alpha} (\hat{x})
:=\hat{\psi}_{\alpha_1}^{\alpha_1}(\hat{x}_1)
\hat{\psi}_{\alpha_2}^{\alpha_2}(\hat{x}_2) \cdots
\hat{\psi}_{\alpha_d}^{\alpha_d}(\hat{x}_d),
$$
and
$$
\hat{\psi}_{\alpha_i}^{\alpha_i}(\hat{x}_i)
:=\frac{1}{2^{\alpha_i} \alpha_i!}
\frac{d^{\alpha_i}}{d \hat{x}_i^{\alpha_i}} (\hat{x}_i+1)^{\alpha_i}
(\hat{x}_i-1)^{\alpha_i}
$$
denotes the $\alpha_i$-th one-dimensional Legendre polynomial of degree
$\alpha_i$ w.r.t.\ $\hat{x}_i\in I,$ $1\le i\le d.$
In this case and with an appropriate indexing, we see that
$P_h^k v=\sum_{i=1}^{\ndofk}(\rho_i^J)^{-1} (v,\psi_i)_{0,T} \psi_i.$

Summarized representations about Legendre polynomials can be found in
\cite[App.\ A]{Karniadakis:05} or \cite[Ch.\ 4]{Quarteroni:94}.

From the definition of the $L^2$-projection, the following properties
easily follow:
\begin{eqnarray}
P_h^k v&=&v \quad \forall v \in \mathbb{Q}_k(T),
\label{L2_proj_eig_inv} \\
\|P_h^k v\|_{0,2,T} &\le & \|v\|_{0,2,T}\quad \forall v \in L^2(T).
\nonumber %
\end{eqnarray}
\begin{lemma}\label{projektionsfehler} %
Let $T$ be an element of an affine partition $\mathcal{T}_h$
such that the corresponding family of partitions $\{\mathcal{T}_h\}_{h>0}$
is locally quasiuniform.
Then, for $1\le l \le k+1,$ $l \in \N,$
there exist constants $C > 0$ independent of $h_T$ and $k$ such that
\begin{eqnarray}
|v - P_h^k v|_{r,2,T}
&\le& C \frac{h_T^{l-r}}{k^{e(r,l)}} |v|_{l,2,T} ,\quad 0\le r \le l,
\nonumber\\
|v - P_h^k v|_{r,2,E}
&\le& C \frac{h_T^{l-1/2-r}}{k^{e(r+1/2+\eps,l)}}|v|_{l,2,T},
\quad 0\le r < l-\frac{1}{2},
\label{L2_proj_error_rand}
\end{eqnarray}
for all $v \in W^{l,2}(T),$
where $\eps>0 $ is arbitrarily small and
$$
e(r,l)=:
\left\{
\begin{array}{ll}
l+1/2-2r, & r \ge 1, \\
l-3r/2, & 0 \le r \le 1. \\
\end{array}
\right.
$$
\end{lemma}
\proof
By \cite[Thm.\ 2.4]{Canuto:82}, there exists a constant $C>0$ independent of $k$
such that, for $r,l \in \R_+,$ $r \le l,$
$$
\|\hat{v}-\hat{P}_h^k\hat{v} \|_{r,2,\hat{T}} \le C k^{- e(r,l)}
\|\hat{v}\|_{l,2,\hat{T}} \quad \forall \hat{v}\in W^{l,2}(\hat{T}).
$$
As in the proof of Lemma \ref{interpolationsfehler},
this estimate together with (\ref{L2_proj_eig_inv})
and the Deny-Lions lemma (\cite[Thm.\ 14.1]{Ciarlet:91a})
implies that, for $l \le k+1,$
$$
|\hat{v}-\hat{P}_h^k\hat{v}|_{r,2,\hat{T}} \le C k^{-e(r,l)}|\hat{v}|_{l,2,\hat{T}}
\quad \forall \hat{v}\in W^{l,2}(\hat{T}),
$$
where $C=C(d,l,\hat{T})>0.$
In particular, the constant does not depend on the polynomial degree $k.$

At the faces of $\hat{T}$ we make use of the following argument.
Let $\hat{E}\subset\partial\hat{T}$ be an arbitrary but fixed face.
Without loss of generality we may consider it as a subset of $\R^{d-1},$
where the elements $\hat{x}\in\R^{d-1}$ are characterized by $x_d=0$
and the elements of $\hat{T}$ satisfy the condition $x_d\ge 0$
(otherwise we apply a rotation and a translation of the coordinate system;
both operations do not affect the differentiability properties of
the elements of the function spaces under consideration).

Now, let $\hat{w}\in W^{s+1/2,2}(\hat{T})$ for some $s\in(0,l-1/2].$
Then, by Kalyabin's results on Stein's extension operator
(\cite[p.\ 185 and Thm.\ 1]{Kalyabin:85}),
there exists a total extension operator
$E_{\hat{T}}:\;W^{s+1/2,2}(\hat{T})\to W^{s+1/2,2}(\R^d)$
such that
$$
\|E_{\hat{T}}\hat{w}\|_{s+1/2,2,\R^d}\le C \|\hat{w}\|_{s+1/2,2,\hat{T}}.
$$
The trace theorem (\cite[Thm.\ 7.43 together with Rem.\ 7.33]{Adams:03})
implies that
$$
\|E_{\hat{T}}\hat{w}\|_{s,2,\R^{d-1}}\le C \|E_{\hat{T}}\hat{w}\|_{s+1/2,2,\R^d},
$$
so
$$
\|E_{\hat{T}}\hat{w}\|_{s,2,\hat{E}}
\le \|E_{\hat{T}}\hat{w}\|_{s,2,\R^{d-1}}
\le C \|\hat{w}\|_{s+1/2,2,\hat{T}}.
$$
Note that for small $s\in(0,1),$ we have by a direct trace theorem
for Lipschitz domains (\cite[Thm.\ 3.1]{Jerison:95}) that
$$
\|\hat{w}\|_{s,2,\hat{E}}
\le \|\hat{w}\|_{s,2,\partial\hat{T}}
\le C \|\hat{w}\|_{s+1/2,2,\hat{T}}.
$$
Therefore, $E_{\hat{T}}\hat{w}|_{\hat{E}}=\hat{w}|_{\hat{E}}$
in the sense of $W^{s,2}(\hat{E})$
(consider, for $l\ge 2,$ $\hat{w}\in W^{s+1/2,2}(\hat{T})$ for $s\in[1,l-1/2]$
as an element of, say, $\hat{w}\in W^{1,2}(\hat{T})$)
and we finally arrive at the estimate
\begin{equation}\label{eq:reftracethm}
\|\hat{w}\|_{s,2,\hat{E}}
\le \|E_{\hat{T}}\hat{w}\|_{s,2,\R^{d-1}}
\le C \|\hat{w}\|_{s+1/2,2,\hat{T}},
\quad
s\in(0,l-1/2].
\end{equation}
In summary, for $r\in[0,l-1/2)$ and $\eps>0 $ arbitrarily small, we get
$$
\|\hat{w}\|_{r,2,\hat{E}}
\le \|\hat{w}\|_{r+\eps,2,\hat{E}}
\le C \|\hat{w}\|_{r+\eps+1/2,2,\hat{T}}.
$$
In particular, for $\hat{w}:=\hat{v}-\hat{P}_h^k\hat{v},$ we have
$$
\|\hat{v}-\hat{P}_h^k\hat{v}\|_{r,2,\hat{E}}
\le C \|\hat{v}-\hat{P}_h^k\hat{v}\|_{r+\eps+1/2,2,\hat{T}},
$$
and, by \cite[Thm.\ 2.4]{Canuto:82},
$$
\|\hat{v}-\hat{P}_h^k\hat{v}\|_{r,2,\hat{E}}
\le C k^{-e(r+\eps+1/2,l)}\|\hat{v}\|_{l,2,\hat{T}}.
$$
Thus, by (\ref{L2_proj_eig_inv})
and the Deny-Lions lemma (\cite[Thm.\ 14.1]{Ciarlet:91a}),
$$
|\hat{v}-\hat{P}_h^k\hat{v}|_{r,2,\hat{E}}
\le C k^{-e(r+\eps+1/2,l)}|\hat{v}|_{l,2,\hat{T}}.
$$
The back-transformation to $T$ runs analogously as in the proof of
Lemma \ref{interpolationsfehler}.
\close

In contrast to the estimates of the projection error presented above,
the next assertion can be proved without the use of the Deny-Lions lemma.
\begin{lemma}\label{l:projerr2} %
Let $T$ be an element of an affine partition $\mathcal{T}_h$
such that the corresponding family of partitions $\{\mathcal{T}_h\}_{h>0}$
is locally quasiuniform.
Then, for $0 \le l$ and $v \in W^{l,2}(T)$ with $1\le s\le\min\{k+1,l\},$
there exist constants $C>0$ independent of $h_T$ and $k$ such that
\begin{eqnarray*}
|v -P_h^k v|_{0,2,T} &\le & C \left(\frac{h_T}{k}\right)^{s} |v|_{s,2,T},\\
|v - P_h^k v|_{1,2,T} &\le & C \frac{h_T^{ s -1}}{k^{s-3/2}} |v|_{s,2,T},\\
|v-P_h^k v|_{0,2,\partial T}&\le & C \left(\frac{h_T}{k}\right)^{s-1/2} |v|_{s,2,T}.
\end{eqnarray*}
\end{lemma}
\proof
\cite[Cor.\ 3.15, (3.1)-(3.3)]{Georgoulis:03},
\cite[Lemma 3.5]{Houston:02},
\cite[Lemma 6.1, Rem.\ 6.2]{Sudirham:06}.
\close

\begin{remark} %
For comparison, we mention the following special case of
Lemma \ref{projektionsfehler}, (\ref{L2_proj_error_rand}):
$$
|v - P_h^k v|_{0,2,E} \le C \frac{h_T^{l-1/2}}{k^{l-3/2(1/2+\eps)}}
|v|_{l,2,T}
\stackrel{\eps < 1/6}{\le }
C \frac{h_T^{l-1/2}}{k^{l-1}}|v|_{l,2,T}\,.
$$
We see that this estimate is suboptimal of order $k^{1/2}.$
The reason lies in the behavior of the $L^2$-projection error
measured in the $|\cdot|_{r,2,T}$-norm for $r>0.$
Indeed, any proof of such an estimate which is based on a trace inequality
will result in right-hand side bounds where the occuring norms
cause a suboptimal result.
To overcome this problem, a direct estimate of the $L^2$-projection error
is helpful, see \cite[Section 3.3]{Houston:02}.
\end{remark}

\section{Quadrature and lumping}
\label{sec:quadrature}
In order to evaluate the occuring integrals approximately
we will consider interpolatory quadrature rules (cf. \cite[Def.\ 8.1]{Ern:04}).

\begin{definition}\label{nodal_quad} %
Let $T\subset\R^d$ be a nonempty, compact, connected subset with
a Lipschitzian boundary (cf.\ \cite[Def.\ 4.9]{Adams:03}).
A quadrature rule on $T$ with $\ndofk$ nodes is characterized
\begin{enumerate}
\item by a set consisting of $\ndofk$ real numbers
$\left\{\omega_1^J,\dots,\omega_{\ndofk}^J \right\}$ called
\emph{weights,} and
\item by a set $\mathcal{Q}$ consisting of $\ndofk$ points
$\{x_1,\dots,x_{\ndofk}\}\subset T,$ where $x_i \ne x_j$ if $i \ne j,$
called \emph{quadrature nodes.}
\end{enumerate}
The largest natural number $k$ such that
$$
\int_T p(x) \, dx= \sum_{i=1}^{\ndofk} \omega_i^J p(x_i)
\quad
\forall p \in \mathbb{Q}_k(T)
$$
is called the \emph{degree of precision} or the \emph{quadrature order}
of the quadrature rule.
\end{definition}
From
\begin{eqnarray*}
\int_T p(x) \, dx
&=&\int_{\hat{T}} p(F_T(\hat{x})) |\det J_T| \, d\hat{x}\\
&=& \int_{\hat{T}} \sum_{\alpha \in \N_0^d, \, \|\alpha\|_{\ell^{\infty}} \le k}
p(F_T(\hat{x}_\alpha)) \hat{\varphi}_{\alpha}(\hat{x}) |\det J_T| \, d\hat{x} \\
&=& \sum_{\alpha \in \N_0^d, \, \|\alpha\|_{\ell^{\infty}} \le k} \int_{\hat{T}}
\hat{\varphi}_{\alpha}(\hat{x}) |\det J_T| \, d\hat{x} \, p(x_\alpha)
= \sum_{\alpha \in \N_0^d, \, \|\alpha\|_{\ell^{\infty}} \le k}
\omega_{\alpha}^J \, p(x_\alpha), \quad x_{\alpha} \in \mathcal{Q},
\end{eqnarray*}
we immediately get the weights
$$
\omega_i^J=\int_{\hat{T}} \hat{\varphi}_i (\hat{x}) |\det J_T| \, d\hat{x} =\int_T
\varphi_i(x) \, dx,\quad 1 \le i \le \ndofk,
$$
corresponding to the nodes $x_i:=F_T(\hat{x}_i),$
where $\hat{\varphi}_i$ is the Lagrange basis function to $\hat{x}_i.$
For reasons of numerical robustness of the quadrature rule the following condition
has to be satisfied:
$$
\omega_i^J >0 \quad\mbox{for}\quad 1 \le i \le \ndofk\,.
$$
The Definition \ref{nodal_quad} of the nodal quadrature rule can be used to
define a discrete inner product. The induced discrete norm allows
to derive norm equivalence estimates w.r.t.\ the $L^p$-norm for discrete arguments
such that the equivalence constants depend only on $p,$ the polynomial degree $k$
and the distribution of qudrature nodes.

To do so, we start with the definition of control volumina
and associated lumping operators. The discrete $L^2$-norm is given by
$$
\begin{array}{rcl}
\|v\|_{0,2,T,h}^2
&:=&\displaystyle
(v,v)_{0,T,h}:=\sum_{\alpha \in \N_0^d, \,
\|\alpha\|_{\ell^{\infty}} \le k} \omega_{\alpha}^J \, v(x_\alpha)^2 \\
&=&\displaystyle
\sum_{\alpha_1=0}^k \cdots \sum_{\alpha_d=0}^k \omega_{\alpha_1}^J \cdots
\omega_{\alpha_d}^J v(x_\alpha)^2,\quad x_{\alpha} \in \mathcal{Q},
\end{array}
$$
with $\sum_{\alpha \in \N_0^d, \, \|\alpha\|_{\ell^{\infty}} \le k}
\omega_{\alpha}^J=|T|,$ $\omega_{\alpha} >0,$ $\forall \alpha \in \N_0^d.$
The generalization to the $L^p$-norm for $p\in[1,\infty)$ is straightforward.

The control volumina are introduced as follows:
$$
\Omega_\alpha=\left( \sum_{i=0}^{\alpha_1-1} \omega_{i}^J, \sum_{i=0}^{\alpha_1}
\omega_{i}^J \right) \times \cdots \times \left(\sum_{i=0}^{\alpha_d-1}
\omega_{i}^J, \sum_{i=0}^{\alpha_d} \omega_{i}^J \right),
$$
where we use the convention
$\sum_{i=0}^{-1}\omega_i^J:=0.$
As a consequence, for the $d$-dimensional Jordan measure of $\Omega_{\alpha}$
we have that
\begin{equation}\label{om_a}
|\Omega_{\alpha}|=\omega_{\alpha_1}^J \omega_{\alpha_2}^J \cdots
\omega_{\alpha_d}^J=\omega_{\alpha}^J,
\end{equation}
and
\begin{equation}\label{sum_om}
\sum_{\alpha \in \N_0^d, \, \|\alpha\|_{\ell^{\infty}} \le k}
|\Omega_{\alpha}|=|T|.
\end{equation}
The lumping operator
$L_h:\;C(T)\to L^\infty(T)$ is defined by
$$
L_h v =\sum_{i=1}^{\ndofk} v(x_i) \chi_{\Omega_i}, \quad x_i \in \mathcal{Q},
$$
where
$$
\chi_G(x)=
\left\{\begin{array}{ll}
0, &\quad x \notin G,\\
1, &\quad x \in G
\end{array}\right.
$$
is the indicator function of the set $G.$
Integrating the $p$-th power of the lumping operator, from
(\ref{om_a}), (\ref{sum_om}) together with the affine transformation
of the reference element we see that
\begin{eqnarray}
\|L_h(v)\|_{0,p,T}^p &=& |T| \int_{\hat{T}} |L_h(v)|^p \, d\hat{x}=
\sum_{j=1}^{\ndofk} |\Omega_j|\int_{\hat{T}} |L_h(v)|^p \, d\hat{x}\nonumber\\
&=&\sum_{j=1}^{\ndofk} \int_{\Omega_j} \left| \sum_{i=0}^{\ndofk} v(x_i)
\chi_{\Omega_i}(x) \right|^p \,dx
=\sum_{j=1}^{\ndofk}|\Omega_j| |v(x_j)|^p \label{id_lumping_quad}\\
&=&\|v\|_{0,p,T,h}^p,\quad x_i \in \mathcal{Q}.\nonumber
\end{eqnarray}
Furthermore,
$$
\|v\|_{0,p,T,h}^p
\stackrel{\mathcal{N}=\mathcal{Q}}{=}\int_{T} \sum_{j=1}^{\ndofk} |v(x_j)|^p
\varphi_j(x) \, dx = \int_T I_h^k(|v|^p) \, dx.
$$
\begin{lemma}\label{normaeq_lumping} %
There exist constants $C_{L_1},$ $C_{L_2}>0$ independent of $h_T$
such that the following equivalence estimates are valid:
\begin{equation}\label{normaeq_lumping_eq}
C_{L_1} \|v\|_{0,p,T} \le \|L_h(v)\|_{0,p,T} \le C_{L_2} \|v\|_{0,p,T}
\quad\forall v \in \mathbb{Q}_{k}(T),
\quad
p\in[2,\infty).
\end{equation}
\end{lemma}
\proof
First we mention that, as a consequence of the substitution rule,
the constants do not depend on $h_T$ for all elements $T$ which result
from an affine transformation of the reference element.
Since $\|L_h(v)\|_{0,p,{T}}^p=\sum_{j=1}^{\ndofk} |\Omega_j| |v(x_j)|^p$
and $\|v\|_{\ell^p}^p=\sum_{j=1}^{\ndofk}|v(x_j)|^p,$
we obtain, by Nikolski's lemma (cf.\ Lemma \ref{l:Nikolski})
and H\"older's inequality for sums,
$$
\begin{array}{rcl}
\|\hat{v}\|_{0,p,\hat{T}}
&\le& (3k^2)^{d\frac{p-2}{2p}}\|\hat{v}\|_{0,2,\hat{T}}
\le (3k^2)^{d\frac{p-2}{2p}}\lambda_{\max}(\hat{M})^{1/2}
\|\hat{v}\|_{\ell^2} \\
&\le& (3k^2(k+1))^{d\frac{p-2}{2p}}\lambda_{\max}(\hat{M})^{1/2}
\|\hat{v}\|_{\ell^p} \\
&\le& (3k^2(k+1))^{d\frac{p-2}{2p}} \lambda_{\max}(\hat{M})^{1/2} \left(
\min\limits_{1 \le j \le \ndofk} |\hat\Omega_j|
\right)^{-1/p}\|\hat{L}_h(\hat{v})\|_{0,p,\hat{T}}
\end{array}
$$
and
$$
\begin{array}{rcl}
\|\hat{L}_h(\hat{v})\|_{0,p,\hat{T}}
&\le& \left(\max\limits_{1 \le j \le \ndofk} |\hat\Omega_j|
\right)^{1/p}\|\hat{v}\|_{\ell^p}
\le \left(\max\limits_{1 \le j \le \ndofk} |\hat\Omega_j|
\right)^{1/p}\|\hat{v}\|_{\ell^2} \\
&\le& \lambda_{\min}(\hat{M})^{-1/2}
\left(\max\limits_{1 \le j \le \ndofk} |\hat\Omega_j|\right)^{1/p}\|\hat{v}\|_{0,2,\hat{T}}\\
&\le& \lambda_{\min}(\hat{M})^{-1/2}|\hat{T}|^{(p-2)/2p}
\left(\max\limits_{1 \le j \le \ndofk} |\hat\Omega_j|\right)^{1/p}\|\hat{v}\|_{0,p,\hat{T}},
\end{array}
$$
where $\Lambda_{\min}(\hat{M}),$ $ \Lambda_{\max}(\hat{M})$
denote the minimal and maximal eigenvalues, resp.,
of the corresponding mass matrix with entries
$\hat{M}_{ij}:=\int_{\hat{T}} \hat{\varphi}_i \hat{\varphi}_j \, d\hat{x}.$
\close

In the approximation of the $L^2$-projection by means of quadrature rules,
the choice of the positions of the quadrature nodes plays an essential role.

On the other hand,
there is also some freedom in the choice of the node set $\mathcal{N}$
in the case of Lagrange basis polynomials
$\hat{\varphi}_i^k$ for $0 \le i \le k.$
In the case of coinciding node sets $\mathcal{N}$ and $\mathcal{Q}$
we see that the definitions of the $L^2$-projection
and of the Lagrange basis polynomials result in
$$
0\stackrel{!}{=}(v-P_h^k v,w)_{0,T,h}
=\sum_{i=1}^{\ndofk} \omega_i^J \left\{v(x_i)- (P_h^k v)(x_i) \right\} w(x_i)
\quad \forall w \in \mathbb{Q}_k(T)\,,
$$
that is
$v(x_i)= (P_h^k v)(x_i)$
for all quadrature nodes.

Thus the representation
$(P_h^k v)(x)=\sum_{j=1}^{\ndofk} (P_h^k v)_j\varphi_j(x)$
of the $L^2$-projection leads to
$(P_h^kv)_i=v(x_i).$
Consequently, for the above approximation of the $L^2$-projection,
we have the implication
$\mathcal{N}=\mathcal{Q} \Rightarrow P_h^k=I_h^k.$

The best accuracy can be reached by using the Gauss quadrature rules
in the following sense.
A $(k+1)^d$-node Gauss quadrature rule yields exact results
for polynomials of maximum degree $2k+1$
(see, e.g., \cite[Prop.\ 8.2]{Ern:04}, \cite[Sect.\ 13,14]{Bernardi:97}
or \cite[pp.\ 448]{Canuto:07}).
The quadrature points w.r.t.\ the domain of integration $I=[-1,1]$
are the zeros of the Legendre polynomial $\hat{\psi}_{k+1}^{k+1}(\hat{x}),$
$\hat{x}\in I,$ of degree $k+1.$
However, the use of Gauss quadrature rules does not lead to
optimal estimates of the quadrature error w.r.t.\ the $W^{l,2}$-norm,
as the following result indicates.
\begin{lemma} %
For all real numbers $r$ and $l,$ $r\le l,$ $l\ge 1,$ there exists
a constant $C >0$ independent of $k$ such that
$$
\|\hat{v} - \hat{I}_{G}^k \hat{v}\|_{r,2,I} \le C k^{-e(r,l)}\|\hat{v}\|_{l,2,I}
\quad \forall \hat{v} \in W^{l,2}(I),
$$
where $e(r,l)$ is the function introduced in Lemma \ref{projektionsfehler}.
\end{lemma}
\proof
See \cite[(13.15)]{Bernardi:97}.
\close

Alternatively, Gauss-Lobatto quadrature rules can be considered.
Here, the quadrature points are the zeros of the polynomial
$$
(x+1)(x-1)(\hat{\psi}_k^k)'(x),\quad x\in I\,.
$$
This quadrature rule is exact for polynomials of maximum degree $2k-1$
(see, e.g., \cite[pp.\ 448]{Canuto:07}).
The boundary points are quadrature nodes.
In contrast to the previous lemma, the following result is valid.
\begin{lemma}\label{gl_error} %
For all real numbers $r$ and $l$ such that $2l > d+r$ and $0 \le r \le 1,$
there exists a constant $C >0$ independent of $k$ such that
\begin{equation}\label{gl_error_eq}
\|\hat{v}-\hat{I}_{GL}^k \hat{v}\|_{r,2,\hat{T}}
\le C k^{r-l}\|\hat{v}\|_{l,2,\hat{T}}
\quad \forall \hat{v} \in W^{l,2}(\hat{T}).
\end{equation}
\end{lemma}
\proof
See \cite[Thm.\ 14.2]{Bernardi:97}.
\close

In the application to quadrature, we have the following result
(cf.\ Lemma \ref{normaeq_lumping}).
\begin{lemma} %
Let $\mathcal{Q}$ be the set of Gauss-Lobatto quadrature nodes.
Then there exists a constant $C>0$ independent of $h_T$ and $k$
such that the following equivalence estimates are valid:
\begin{equation}\label{normaeq_lumping_gl}
\|v\|_{0,2,T} \le \|L_h(v)\|_{0,2,T} \le C \|v\|_{0,2,T}
\quad \forall v \in \mathbb{Q}_{k}(T).
\end{equation}
\end{lemma}
\proof
First we prove the statement for the reference element.
On $\hat{T},$ we have
$$
\|L_h(\hat{v})\|_{0,2,\hat{T}}^2=\|\hat{v}\|_{0,2,\hat{T},GL}^2
$$
by (\ref{id_lumping_quad});
the remaining part is a consequence of \cite[(3.9)]{Canuto:82}.
The affine back-transformation to the original element shows that the constants
are independent of $h_T.$
\close
\begin{remark}
Analogously to Lemma \ref{normaeq_lumping} we can conclude that the
above inequalities (\ref{normaeq_lumping_gl}) can be extended to the
general case $p\in[2,\infty).$
\end{remark}
A further interesting property of Gauss-Lobatto quadrature nodes
is related with the square sum of Lagrange polynomials.
\begin{lemma} %
Let $\hat{\varphi}_{\alpha}(\hat{x})$ be the Lagrange polynomials w.r.t.\ the
Gauss-Lobatto quadrature nodes. Then we have the estimate
$$
\sum_{\alpha \in \N_0^d, \atop \|\alpha\|_{\ell^\infty} \le k}
\hat{\varphi}_{\alpha}(\hat{x})^2 \le 1
\quad \forall \hat{x} \in \hat{T}.
$$
\end{lemma}
\proof
The proof is a consequence of the tensor product representation
together with \cite[\S\,2]{Fejer:32}:
$$
\sum_{\alpha \in \N_0^d, \atop \|\alpha\|_{\ell^\infty} \le k}
\hat{\varphi}_{\alpha}(\hat{x})^2=\underbrace{\sum_{\alpha_1=0}^k
\hat{\varphi}_{\alpha_1}^k (\hat{x}_1)^2}_{\le 1} \cdots
\underbrace{\sum_{\alpha_d=0}^k \hat{\varphi}_{\alpha_d}^k
(\hat{x}_d)^2}_{\le 1} \le 1.
\qquad\close
$$
\section{Projection and interpolation errors w.r.t.\ Gauss-Lobatto
quadrature nodes}
\label{sec:errors}
The Legendre polynomials used in the representation of the discrete
$L^2$-projector possess two important properties:
On the one hand, they form an orthogonal basis of $\mathbb{Q}_k(T),$
and, consequently, the corresponding mass matrix is diagonal.
On the other hand, there exists a hierarchical decomposition
of $\mathbb{Q}_k(T)$ in the following sense (cf.\ \cite[Def.\ 1.18]{Ern:04}).
\begin{definition}[Hierarchical modal basis] %
A family $\{\mathcal{B}_k\}_{k\in\N_0},$ where $\mathcal{B}_k$
denotes a set of polynomials, is called \emph{hierarchical modal basis}
if, for all $k\in\N_0,$ the following properties are satisfied:
\begin{enumerate}
\item
$\mathcal{B}_k$ is a basis of $\mathbb{Q}_k,$
\item
$\mathcal{B}_k\subset\mathcal{B}_{k+1}.$
\end{enumerate}
\end{definition}
If the $L^2$-projection is discretized by means of Lagrange polynomials
for the node sets $\mathcal{N}=\mathcal{Q},$
then the discrete $L^2$-orthogonality is conserved but, unfortunately,
the corresponding Lagrange basis is not a hierarchical modal basis.
This problem can be resolved by introducing the following notion.
By $N:\;\mathcal{B}\to\mathcal{N}$ we denote that bijective function
which assigns a set of polynomials $\mathcal{B}$ to the node set
of the corresponding Lagrange basis.
\begin{definition}[embedded hierarchical nodal basis] %
A family $\{\mathcal{B}_j\}_{1\le j\le \ndofk},$
where $\mathcal{B}_j$ denotes a set of polynomials of maximum degree $k,$
is called \emph{embedded hierarchical nodal basis of degree $K,$}
$K\in\N_0,$ if the following properties are satisfied:
\begin{enumerate}
\item
$\tilde\mathcal{B}_K$ with $N(\tilde\mathcal{B}_K)\subset N(\mathcal{B}_k)$
is a basis of $\mathbb{Q}_K,$
\item
$\mathcal{B}_K\subset\mathcal{B}_k,$
$N(\mathcal{B}_K) = N(\tilde\mathcal{B}_K),$
\item
$\mathcal{B}_k$ is a basis of $\mathbb{Q}_k.$
\end{enumerate}
\end{definition}
\begin{example}[$d=1$]
\rule{0ex}{.1ex}

(i) The Lagrange polynomials w.r.t.\ the Gauss-Lobatto quadrature nodes
form an embedded hierarchical nodal basis of degree 1 and 2.
This follows from the fact that
$N(\tilde\mathcal{B}_1):=\{-1,1\}\subset N(\mathcal{B}_k),$ $k\in\N,$ and
$N(\tilde\mathcal{B}_2):=\{-1,0,1\}\subset N(\mathcal{B}_k),$ $k=2,4,6,\ldots$

(ii) The Lagrange polynomials w.r.t.\ the Gauss-Kronrod quadrature nodes
$\{x_i\}_{i=1}^k=:N(\mathcal{B}_k)$
form an embedded hierarchical nodal basis of degree $K$ for $k=2K+2,$ $K\in\N_0.$
Given $K+1$ Gauss quadrature nodes
$\{\tilde{x}_i\}_{i=0}^K=:N(\tilde\mathcal{B}_K),$
Gauss-Kronrod quadrature rules are defined by adding $K+2$ nodes such that
\begin{eqnarray*}
N(\tilde\mathcal{B}_K)&\subset&N(\mathcal{B}_k),\\
\int_I v\,dx&=&\sum_{i=0}^k \omega_i v(x_i)
\qquad\forall v\in\mathbb{Q}_{3K+4}(I)
\end{eqnarray*}
(see, e.g., \cite{Calvetti:00}).
\end{example}
\begin{figure}[ht]
\includegraphics[width=.45\textwidth]{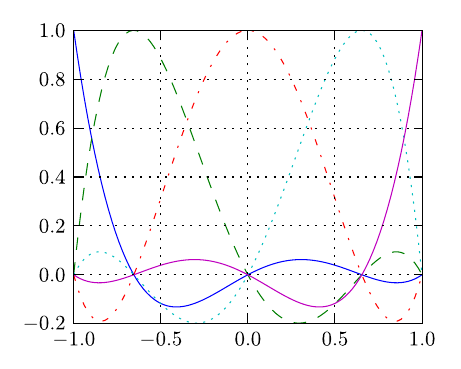}
\quad
\includegraphics[width=.45\textwidth]{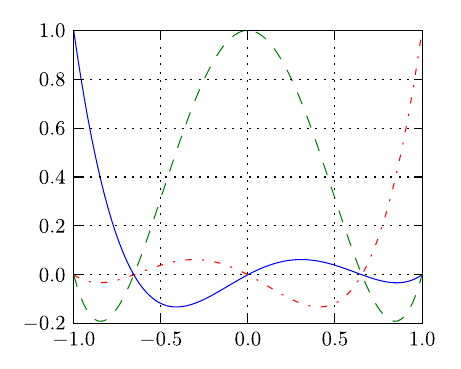}
\caption{Lagrange basis polynomials of $\mathbb{Q}_4(I)$ for
the Gauss-Lobatto nodes (left), corresponding basis polynomials of $V_h^{2,4}(I)$
(right)}
\end{figure}
As an application we investigate some properties of the so-called
\emph{fluctuation operator,} see the definition below.
These results are important in the numerical analysis
of discontinuous Galerkin methods for conservation laws.
\begin{definition} %
Let $\{\mathcal{B}_j\}_{1\le j\le \ndofk}$ be
an embedded hierarchical nodal basis of degree $K,$ $K\in\N_0.$
For $\mathcal{Q}_k(T)=V_h^{K,k}(T)\oplus V_h'(T),$
where $V_h^{K,k}(T):=\spann\mathcal{B}_K$ and
$V_h'(T):=\spann\{\mathcal{B}_k\setminus\mathcal{B}_K\},$
the projector $P_h^{K,k}:\;L^2(T)\to V_h^{K,k}(T)$
is defined by
\begin{equation}\label{L2_embproj} %
(v -P_h^{K,k}v, w)_{0,T,h}=0\quad \forall w \in V_h^{K,k}(T).
\end{equation}
The operator $P_h'(T):=Id-P_h^{K,k}$ in $L^2(T),$ where $Id$
denotes the identity in $L^2(T),$ is called \emph{fluctuation operator.}
\end{definition}
\begin{lemma}\label{l:PhKkerror} %
Let $\{\mathcal{B}_j\}_{1\le j\le \ndofk}$ be
an embedded hierarchical nodal basis of degree $K,$ $K\in\N_0,$
consisting of Lagrange polynomials and suppose $\mathcal{N}=\mathcal{Q}.$
Then there is a constant $C>0$ independent of $h_T$ such that
\begin{eqnarray*}
|v - P_h^{K,k} v|_{r,2,T}
& \le & C h_T^{l-r}\left[|v|_{l,2,T}+|P_h^K v|_{l,2,T}\right], \quad 0 \le r \le l,\\
|v - P_h^{K,k} v|_{r,2,E}
& \le & C h_T^{l-1/2-r}\left[|v|_{l,2,T}+|P_h^K v|_{l,2,T}\right], \quad 0 \le r< l-1/2
\end{eqnarray*}
for all $v\in W^{l,2}(T).$
\end{lemma}
\proof
Under the above assumptions we have, for all Lagrange basis polynomials
$\varphi_j^{K,k}\in \mathcal{B}_k,$ that
$\varphi_j^{K,k}(x_i) = \delta_{ij},$ $1\le i,j \le\ndofK,$
$x_i\in N(\tilde\mathcal{B}_K) = N (\mathcal{B}_K ).$
The definition of the projector $P_h^{K,k}$ and the property
\begin{equation}\label{eq:fop1} %
w \in V_h^{K,k}(T)\ \Longrightarrow\ w(x_i) = 0,
\quad x_i\in N(\mathcal{B}_k)\setminus N(\mathcal{B}_K )
\end{equation}
imply that
\begin{eqnarray*}
0&\stackrel{!}{=}&(v -P_h^{K,k}v, w)_{0,T,h}
=\sum_{i=1}^{\ndofk} \omega_i\left\{v(x_i)-(P_h^{K,k}v)(x_i)\right\} w(x_i)\\
&\stackrel{(\ref{eq:fop1})}{=}&\sum_{i=1}^{\ndofK}
\omega_i\left\{v(x_i)-(P_h^{K,k}v)(x_i)\right\} w(x_i)
\end{eqnarray*}
and, thus,
\begin{equation}\label{eq:fop2} %
v(x_i)\stackrel{!}{=}
\sum_{j=1}^{\ndofK} (P_h^{K,k}v)_j\varphi_j^{K,k}(x_i)
\stackrel{\mathcal{N}=\mathcal{Q}}{=} (P_h^{K,k}v)_i,\quad 1 \le i \le \ndofK\,.
\end{equation}
As a consequence, if $I_h^K$ denotes the Lagrange interpolation operator
w.r.t.\ the nodes $N(\mathcal{B}_K),$
we obtain the relation
$$
I_h^K(P_h^{K,k} v)
=\sum_{i=1}^{\ndofK}(P_h^{K,k}v)(x_i)\varphi_i^K
=\sum_{i,j=1}^{\ndofK}(P_h^{K,k}v)_j\varphi_j^{K,k}(x_i)\varphi_i^K
=\sum_{i=1}^{\ndofK}(P_h^{K,k}v)_i\varphi_i^K
\stackrel{(\ref{eq:fop2})}{=}I_h^K(v)\,.
$$
It can be used, together with Lemma \ref{interpolationsfehler},
in the estimation of the projection error as follows:
\begin{eqnarray*}
|v - P_h^{K,k} v|_{r,2,T}
&\le& |v - I_h^K v|_{r,2,T}+ |I_h^K v - P_h^{K,k} v|_{r,2,T}\\
&\le& |v - I_h^K v|_{r,2,T}+ |P_h^{K,k} v - I_h^K(P_h^{K,k}) v|_{r,2,T}\,.
\quad\close
\end{eqnarray*}
\begin{lemma} %
Let $\{\mathcal{B}_j\}_{1\le j\le \ndofk}$ be
an embedded hierarchical nodal basis of degree $K,$ $K\in\N_0,$
consisting of Lagrange polynomials and suppose $\mathcal{N}=\mathcal{Q}.$
Then:
\begin{rmnumerate}
\item
$P_h^{K,k}$ is a linear operator,
\item
$\|P_h^{K,k}v\|\le C(K,k,\|\cdot\|)\|v\|_{0,\infty,T}$ for an arbitrary norm
$\|\cdot\|$ on $V_h^{K,k}(T).$
\item
For $m\in\N,$ $v\in\mathbb{Q}_k(T),$ the following estimate is valid:
$$
(P_h'(\nabla v),P_h'(\nabla v^{2m-1}))_{0,T,h}\ge 0.
$$
\end{rmnumerate}
\end{lemma}
\proof
The assertion (i) easily follows from (\ref{eq:fop2}):
$$
P_h^{K,k}(\alpha v + \beta w)
=\sum_{i=1}^{\ndofK}\left[\alpha v(x_i) + \beta w(x_i)\right]\varphi_i^{K,k}
=\alpha P_h^{K,k}v + \beta P_h^{K,k}w\,.
$$
Item (ii) results from the estimates
\begin{eqnarray*}
\|P_h^{K,k} v\|
&\le&\sum_{i=1}^{\ndofK}\|v(x_i)\varphi_i^{K,k}\|
\le \sum_{i=1}^{\ndofK}|v(x_i)|\|\varphi_i^{K,k}\|\\
&\le&\left(\sum_{i=1}^{\ndofK}\|\varphi_i^{K,k}\|\right)\|v\|_{0,\infty,T}
\le C(K,k,\|\cdot\|)\|v\|_{0,\infty,T}\,.
\end{eqnarray*}
To verify (iii), we mention the following elementary inequality on $T$ ($p:=2m$):
$$
\nabla v\cdot\nabla v^{p-1}
=\frac{4(p-1)}{p^2}\nabla v^{p/2}\cdot\nabla v^{p/2}\ge 0\,.
$$
Then
\begin{eqnarray*}
(P_h'(\nabla v), P_h'(\nabla v^{p-1}))_{0,T,h}
&\stackrel{(\ref{L2_embproj})}{=}&
(\nabla v, \nabla v^{p-1})_{0,T,h}
- (P_h^{K,k} (\nabla v), P_h^{K,k}(\nabla v^{p-1}))_{0,T,h}\\
&\stackrel{(\ref{eq:fop2})}{=}&
\sum_{i=1}^{\ndofk} \omega_i^J \nabla v(x_i)\cdot\nabla v^{p-1}(x_i)
-\sum_{i=1}^{\ndofK} \omega_i^J \nabla v(x_i)\cdot\nabla v^{p-1}(x_i)\\
&=&\sum_{i=\ndofK+1}^{\ndofk} \omega_i^J \nabla v(x_i)\cdot\nabla v^{p-1}(x_i)\\
&=&\frac{4(p-1)}{p^2}\sum_{i=\ndofK+1}^{\ndofk}
\omega_i^J\nabla v^{p/2}(x_i)\cdot\nabla v^{p/2}(x_i)\ge 0\,.
\end{eqnarray*}
Interchanging the roles of $P_h'$ and $P_h^{K,k},$ we obtain analogously
$(P_h^{K,k}(\nabla v),P_h^{K,k}(\nabla v^{2m-1}))_{0,T,h}\ge 0.$
\close

Motivated by Lemma \ref{l:PhKkerror} and the error estimate (\ref{gl_error_eq})
we investigate the local interpolation error of the Gauss-Lobatto
interpolation operator.
\begin{lemma}\label{l:IGLkerror} %
Let $T$ be an element of an affine partition $\mathcal{T}_h$
such that the corresponding family of partitions $\{\mathcal{T}_h\}_{h>0}$
is locally quasiuniform.
Assume that
$1\le l \le k+1,$ $l \in \N,$ $0 \le r \le 1$ and $2l > d+r.$
Then, for the Lagrange interpolation operator $I_{GL}^k,$
there exist constants $C > 0$ independent of $h_T$ and $k$ such that
\begin{eqnarray}
|v - I_{GL}^k v|_{r,2,T} &
\le & C \left(\frac{h_T}{k}\right)^{l-r} |v|_{l,2,T}, \label{eq:IGLkerror1}\\
|v - I_{GL}^k v|_{r,2,E} &
\le & C \left(\frac{h_T}{k}\right)^{l-1/2-r} |v|_{l,2,T} \label{eq:IGLkerror2}
\end{eqnarray}
for all $v \in W^{l,2}(T).$
\end{lemma}
\proof
Using Lemma \ref{gl_error}, the proof of (\ref{eq:IGLkerror1}) runs analogously
to the proofs of Lemmata \ref{interpolationsfehler} and \ref{projektionsfehler}.

To prove (\ref{eq:IGLkerror2}), we make use of the fact that,
for tensor product elements, the restriction of the Gauss-Lobatto
interpolant to a face $E$ of $T$ is a $(d-1)$-dimensional Gauss-Lobatto
interpolant. From Lemma \ref{gl_error} we see that
$$
|\hat{v}-\hat{I}_{GL}^k \hat{v}|_{r,2,\hat{E}}
\le C k^{r-\tilde{l}}\|\hat{v}\|_{l,2,\hat{E}}
$$
for all $\hat{v} \in W^{\tilde{l},2}(\hat{E})$
and all real numbers $r$ and $\tilde{l}$ such that
$0 \le r \le 1,$ $2\tilde{l}> d-1+r.$

The trace theorem (cf.\ (\ref{eq:reftracethm}) in the proof
of Lemma \ref{projektionsfehler}) leads to
$$
\|\hat{v}\|_{\tilde{l},2,\hat{E}}
\le C \|\hat{v}\|_{\tilde{l}+1/2,2,\hat{T}}.
$$
If in addition $l=\tilde{l}+1/2\in\N,$
the Deny-Lions lemma (\cite[Thm.\ 14.1]{Ciarlet:91a}) implies that
$$
|\hat{v}-\hat{I}_{GL}^k \hat{v}|_{r,2,\hat{E}}
\le C k^{r-l+1/2}|\hat{v}|_{l,2,\hat{T}}.
$$
Turning back to the original variables, we obtain (\ref{eq:IGLkerror2}).
\close
\begin{corollary} %
Let $1\le l \le K + 1.$
Under the assumptions of Lemmata \ref{l:PhKkerror} and \ref{l:IGLkerror},
there exist constants $C>0$ independent of $h_T$ and $k$ such that
\begin{eqnarray*}
|v - P_{GL}^{K,k} v|_{r,2,T}
& \le & C \left(\frac{h_T}{K}\right)^{l-r}\left[|v|_{l,2,T}+|P_{GL}^K v|_{l,2,T}\right],\\
|v - P_{GL}^{K,k} v|_{r,2,E}
& \le & C \left(\frac{h_T}{K}\right)^{l-1/2-r}
\left[|v|_{l,2,T}+|P_{GL}^K v|_{l,2,T}\right]
\end{eqnarray*}
for all $v\in W^{l,2}(T).$
\end{corollary}
\begin{remark}
As in the proof of Lemma \ref{l:PhKkerror},
the error estimate will be reduced to an estimate of a suitable seminorm.
\end{remark}
In general, the operators $\nabla$ and $I_{GL}^k$do not commute.
An estimate of the commutation error is given in the next result.
\begin{lemma} %
Under the assumptions of Lemma \ref{l:IGLkerror},
there exists a constant $C>0$ independent of $h_T$ and $k$ such that
$$
\|I_{GL}^k(\nabla v)-\nabla I_{GL}^k v\|_{0,2,T,GL}
\le C \left(\frac{h_T}{k}\right)^{l-1}|v|_{l,2,T}.
$$
\end{lemma}
\proof
\begin{eqnarray*}
\|I_{GL}^k(\nabla v)-\nabla I_{GL}^k v\|_{0,2,T,GL}
&\stackrel{(\ref{normaeq_lumping_eq})}{\le }&
C \|I_{GL}^k(\nabla v)-\nabla I_{GL}^k v\|_{0,2,T}\\
&\le& C \|\nabla(v-I_{GL}^k v)\|_{0,2,T}
+ C \|\nabla v - I_{GL}^k(\nabla v)\|_{0,2,T}\\
&\stackrel{(\ref{eq:IGLkerror1})}{\le }&
C\left(\frac{h_T}{k}\right)^{l-1} |v|_{l,2,T}
+C\left(\frac{h_T}{k}\right)^{l-1} |\nabla v|_{l-1,2,T}\\
&\le& C\left(\frac{h_T}{k}\right)^{l-1} |v|_{l,2,T}.
\qquad\close
\end{eqnarray*}

\section{Inverse inequalities}
\label{sec:inverse}
Inverse inequalities are an important tool in the analysis
of finite element methods.
Using particular properties if the tensor product representation,
in this section
we derive sharp inverse inequalities such that, with exception of the difference
of the derivative orders and the spatial dimension, all other information is
known explicitly.
It is well-known that in the Nikolski inequality the dependence
of the degree of polynomials cannot be improved, as the example
\begin{equation}\label{eq:Timanex}
p(x)=\left( \frac{1- T_n^2(x) }{1-x^2}\right)^2
\end{equation}
with $T_n$ being the $n$-th \v{C}eby\v{s}ov polynomial of the first kind
demonstrates, see \cite[p.\ 263]{Timan:63}
(although in other papers, e.g.
\cite[Lemma 2.4]{Gudi:08} or \cite[Lemma 1]{Chen:93b},
different statements can be found).
The representation
$$
T_n(x)=\frac{1}{2} \left[ (x +\sqrt{x^2-1})^n + (x -\sqrt{x^2-1} )^n \right],
$$
shows that $1-T_n^2(x)$ has the zeros $-1$ and $1,$ thus
(\ref{eq:Timanex}) is really a polynomial.

Now we present a generalization of Nikolski and Markov inequalities
(see \cite{Nikolski:51}, \cite[Sect.\ III]{Hille:37})
to the tensor product situation, and we include it into the proof
of an inverse inequality given in \cite[Lemma 1.138]{Ern:04}.

Throughout this section we restrict ourselves to the case $k\ge 1.$

\begin{lemma}[local inverse inequality]\label{inverse_estimate} %
Given a reference element $\{\hat{T},\hat{P},\hat{\Sigma} \}$ such that,
for $l \ge 0,$ the embedding $\hat{P} \subset W^{l,\infty}(\hat{T})$
is satisfied.
Let $\{\mathcal{T}_h\}_{h\in(0,1]}$ be a locally quasiuniform family
of affine partitions. If $0 \le m \le l,$ then there exist
\begin{enumerate}
\item
for $1 \le p,q \le \infty$
a constant $C=C(l,m,p,q,d,\sigma_0,\hat{T},P(\hat{T}))>0$
such that
\begin{equation}\label{inv_allg}
\|v\|_{l,p,T} \le C h_T^{m-l+d(\frac{1}{p}-\frac{1}{q})} \|v\|_{m,q,T}
\qquad \forall v \in P(T),
\end{equation}
\item
for $1 \le q \le p \le \infty$ and
$\hat{P}:=\mathbb{Q}_k(\hat{T})$
a constant $C=C(l,m,p,d,\sigma_0)>0$
such that
\begin{equation}\label{inv_lagrange}
\|v\|_{l,p,T}
\le C \left(\frac{h_T}{k^2}\right)^{m-l}
\left( \frac{h_T}{2 (q+1)k^2} \right)^{d(\frac{1}{p}-\frac{1}{q})} \|v\|_{m,q,T}
\qquad \forall v \in P(T).
\end{equation}
\end{enumerate}
\end{lemma}
The proof will be given after the presentation of the above mentioned
Nikolski and Markov inequalities.
\begin{lemma}[Nikolski]\label{l:Nikolski} %
For $0 < q \le p \le \infty$ and $\hat{v} \in \mathbb{Q}_k(\hat{T}),$
the following estimate is valid:
\begin{equation}\label{nikolskii_ungleichung}
\|\hat{v}\|_{0,p,\hat{T}}
\le \left( (q+1)k^2 \right)^{-d(1/p -1/q)} \|\hat{v}\|_{0,q,\hat{T}}.
\end{equation}
\end{lemma}
\proof
Due to \cite[Thm.\ 4.2.6]{DeVore:93}, in the one-dimensional situation
we have the following inequality:
$$
\|\hat{v}\|_{0,p,I}
\le \left((q+1)k^2\right)^{-(1/p -1/q)}\|\hat{v}\|_{0,q,I}
\quad\forall \hat{v}\in\mathbb{Q}_k(I),
$$
where, as above, $I=[-1,1].$
Then, for $1 \le i \le d,$ we also have that
$$
\begin{array}{rcl}
&&\|\hat{v}(\hat{x}_1,\dots,\hat{x}_{i-1},\cdot,
\hat{x}_{i+1},\dots,\hat{x}_d)\|_{0,\infty,I} \\
&\le& \left( (q+1)k^2 \right)^{1/q}
\|\hat{v}(\hat{x}_1,\dots,\hat{x}_{i-1},\cdot,
\hat{x}_{i+1},\dots,\hat{x}_d)\|_{0,q,I}.
\end{array}
$$
A successive application of this inequality leads to
$$
\begin{array}{rcl}
\|\hat{v}\|_{0,\infty,\hat{T}}^q
&=& \max\limits_{\hat{x}_1 \in I} \cdots \max\limits_{\hat{x}_d \in I}
|\hat{v} (\hat{x}_1,\dots,\hat{x}_d)|^q \\
&\le& \max\limits_{\hat{x}_1 \in I} \cdots \max\limits_{\hat{x}_{d-1} \in I}
((q+1)k^2) \int_I |\hat{v} (\hat{x}_1,\dots,\hat{x}_d)|^q d\hat{x}_d \\
&\le& ((q+1)k^2) \max\limits_{\hat{x}_1 \in I} \cdots \max\limits_{\hat{x}_{d-2} \in I}
\int_I \max\limits_{\hat{x}_{d-1}} |\hat{v} (\hat{x}_1,\dots,
\hat{x}_d)|^q \,d\hat{x}_d \\
&\le& ((q+1)k^2)^2 \max\limits_{\hat{x}_1 \in I}
\cdots \max\limits_{\hat{x}_{d-2} \in I}
\int_I \int_I |\hat{v} (\hat{x}_1,\dots,
\hat{x}_d)|^q \,d\hat{x}_{d-1} d\hat{x}_d \\
&\le& \cdots \le ((q+1)k^2)^d \|\hat{v} \|_{0,q,\hat{T}}^q,\quad 0 <q<\infty. \\
\end{array}
$$
It remains to make use of an $L^p$-type interpolation inequality
(see, e.g., \cite[Cor.\ B.7]{Ern:04}):
$$
\begin{array}{rcl}
\|\hat{v}\|_{0,p,\hat{T}}
&\le& \|\hat{v}\|_{0,q,\hat{T}}^{\frac{q}{p}}
\|\hat{v}\|_{0,\infty,\hat{T}}^{1-\frac{q}{p}} \\
&\le& \|\hat{v}\|_{0,q,\hat{T}}^{\frac{q}{p}}
\left( (q+1)k^2 \right)^{\frac{d}{q}(1-\frac{q}{p})}
\|\hat{v}\|_{0,q,\hat{T}}^{1-\frac{q}{p}}\\
&\le& \left( (q+1)k^2 \right)^{-d(\frac{1}{p}-\frac{1}{q} )}\|\hat{v}\|_{0,q,\hat{T}},
\quad 0 < q \le p \le \infty.
\qquad\close
\end{array}
$$
\begin{lemma}[generalized Markov inequality]\label{l:genMarkov} %
For $v \in \mathbb{Q}_k(I)$ and $p\ge 1,$ the following estimate is valid:
$$
\|v'\|_{0,p,I} \le C_M(p) k^2 \|v\|_{0,p,I},
$$
where
$$
C_{M,p}=C_M(p):=2 (p-1)^{1/p-1}
\left( p + \frac{1}{k} \right)\left( 1+\frac{p}{kp-p+1} \right)^{k-1+1/p}.
$$
For $p=\infty,$ we even have that
$\displaystyle
C_{M,\infty}=1 < \lim_{p\to \infty} C_M(p)=2e.$
Furthermore, $\forall p \in \N: C_M(p) \le C_M:=6 e^{1+1/e}.$
\end{lemma}
\proof
For $1 < p < \infty,$ the assertion is proved in \cite[Sect.\ III]{Hille:37}.
The estimate $C_M(p) \le C_M=6 e^{1+1/e}$ for all $p\in\N$
can be found in \cite[p.\ 590]{Milovanovic:94}.
For $p=\infty,$ the result follows from \cite[Thm.\ 4.1.4]{DeVore:93}.
\close

\begin{remark}
The constants $C_M(p)$ can be improved.
For instance, in \cite[Cor.\ 2.10]{Baran:98} the existence of a constant
$\tilde{C}_M(p)$ such that
$\lim\limits_{p \to \infty} \tilde{C}_M(p)=1,$ $p > 2$
is demonstrated.
\end{remark}

\begin{corollary}
For $1 \le p \le \infty$ and $\hat{v} \in \mathbb{Q}_k(\hat{T}),$
the following estimate is valid:
\begin{equation}\label{inv_l_0}
\|\hat{v}\|_{l,p,\hat{T}} \le {d+l \choose l}^{\frac{1}{p}}
\left( C_{M,p}\,k^2 \right)^l \|\hat{v}\|_{0,p,\hat{T}}
\end{equation}
(with the constant $C_{M,p}$ as in Lemma \ref{l:genMarkov}).
\end{corollary}
\proof
Iterating over the order of derivatives, it is sufficient to verify
the estimate
$$
\|\partial^\alpha \hat{v}\|_{0,p,\hat{T}}
\le ( C_{M,p} k^2 )^{|\alpha|} \|\hat{v}\|_{0,p,\hat{T}}
$$
for $|\alpha|=1.$
To do so, set
$\alpha :=e_i,$ where $e_i$ denotes the $i$-th coordinate unit vector.

For $p < \infty$ and $1 \le i \le d$ the generalized Markov inequality
implies that
$$
\begin{array}{rcl}
&&\|\hat{v}'(\hat{x}_1,\dots,\hat{x}_{i-1},\cdot,\hat{x}_{i+1},\dots,\hat{x}_d)
\|_{0,p,I}^p \\
&\le& C_{M,p}^p k^{2p}
\|\hat{v}(\hat{x}_1,\dots,\hat{x}_{i-1},
\cdot,\hat{x}_{i+1},\dots,\hat{x}_d)\|_{0,p,I}^p,
\end{array}
$$
hence, on $\hat{T},$
$$
\begin{array}{rcl}
\|\partial^\alpha \hat{v} \|_{0,p,\hat{T}}^p
&=& \int_I \dots \int_I |\partial_i \hat{v} (\hat{x})|^p
\, d\hat{x}_i d\hat{x}_1 \dots d\hat{x}_{i-1} d\hat{x}_{i+1} \dots d\hat{x}_d\\
&\le& C_{M,p}^p k^{2p} \|\hat{v}\|_{0,p,\hat{T}}^p.
\end{array}
$$
For $p=\infty,$ there is a point $\hat{y}\in\hat{T}$ such that
$$
\begin{array}{rcl}
\|\partial^\alpha \hat{v} \|_{0,\infty,\hat{T}}
&=&|\partial_i \hat{v}(\hat{y})|
= \max\limits_{\hat{x}_i \in I}
|\partial_i \hat{v}(\hat{y}_1,\dots,\hat{y}_{i-1},\hat{x}_i,\hat{y}_{i+1},\dots,
\hat{y}_d)| \\
&\le& C_{M,\infty} k^2
\|\hat{v}(\hat{y}_1,\dots,\hat{y}_{i-1},\cdot,\hat{y}_{i+1},\dots,
\hat{y}_d)\|_{0,\infty,I} \\
&\le& C_{M,\infty} k^2 \|\hat{v}\|_{0,\infty,\hat{T}}.
\end{array}
$$
Consequently,
$$
\begin{array}{rcl}
\|\hat{v}\|_{l,p,\hat{T}}^p
&=& \sum_{i=0}^l \sum_{|\alpha|=i} \|\partial^\alpha \hat{v}\|_{0,p,\hat{T}}^p
\le \sum_{i=0}^l \sum_{|\alpha|=i}
\left( C_{M,p} k^2 \right)^{|\alpha|p}\|\hat{v}\|_{0,p,\hat{T}}^p\\
&=& \sum_{i=0}^l {d+i-1 \choose i}
\left( C_{M,p} k^2 \right)^{ip} \|\hat{v}\|_{0,p,\hat{T}}^p \\
&\le& {d+l \choose l} \left( C_{M,p} k^2 \right)^{lp} \|\hat{v}\|_{0,p,\hat{T}}^p,
\end{array}
$$
and
$$
\begin{array}{rcl}
\|\hat{v}\|_{l,\infty,\hat{T}}
&=&\max\limits_{0 \le |\alpha| \le l}\|\partial^\alpha \hat{v}\|_{0,\infty,\hat{T}}
\le \max\limits_{0 \le |\alpha| \le l}\left( C_{M,\infty} k^2\right)^{|\alpha|}
\| \hat{v}\|_{0,\infty,\hat{T}} \\
&\le& \left( C_{M,\infty} k^2\right)^{l} \| \hat{v}\|_{0,\infty,\hat{T}}.
\qquad\close
\end{array}
$$
Now we are ready to prove the inverse estimates
(\ref{inv_allg}) and (\ref{inv_lagrange}).

\textbf{Proof} of Lemma \ref{inverse_estimate}:

The first estimate is proved along the lines of the proof of
\cite[Lemma 1.138]{Ern:04}.
In the particular case of a Lagrange reference element we use the
same idea of proof but make use of the results prepared in this section.

As usual, first the assertion is proved for $m=0.$
By $(\ref{inv_l_0})$ and $(\ref{nikolskii_ungleichung}),$ on $\hat{T}$
we have that
\begin{equation}\label{inv_ref}
\|\hat{v}\|_{l,p,\hat{T}} \le {d+l \choose l}^{\frac{1}{p}}
\left( C_{M,p} k^2 \right)^l \left( (q+1)k^2 \right)^{-d(\frac{1}{p}-\frac{1}{q})}
\|\hat{v}\|_{0,q,\hat{T}}.
\end{equation}
In order to get the corresponding estimate for the transformed element $T,$
we make use of
$(\ref{trans_ungl2}),$ $(\ref{help_scaling2})$ and $(\ref{help_scaling3})$ for
$0 \le j \le l:$
$$
\begin{array}{rcl}
| v|_{j,p,T}^p
&\le& \left\{ C_{j,d} \|J_T^{-1}\|_{\ell^2}^{j} |\det J_T|^{1/p}\right\}^p
|\hat{v}|_{j,p,\hat{T}}^p \\
&\le& \left\{ C_{j,d} \left(2 \sqrt{d} \sigma_0 h_T^{-1}\right)^{j}
\left(\frac{h_T}{2}\right)^{\frac{d}{p}} \right\}^p \|\hat{v}\|_{j,p,\hat{T}}^p \\
&\hspace{-2mm}\stackrel{(\ref{inv_ref})}{\le}&
\left\{{d+j \choose j}^{\frac{1}{p}} \tilde{C}(j,p,d,\sigma_0)
\left( \frac{h_T}{k^2} \right)^{-j} \left(\frac{h_T}{2}\right)^{\frac{d}{p}}
\left( (q+1)k^2 \right)^{-d(\frac{1}{p}-\frac{1}{q})} \right\}^p
\|\hat{v}\|_{0,q,\hat{T}}^{p},
\end{array}
$$
where $\tilde{C}(j,p,d,\sigma_0)= C_{j,d} \left(2\sqrt{d} \sigma_0 C_{M,p} \right)^{j}.$
By $(\ref{trans_ungl1})$ with $l=0,$ the back-transformation yields
$$
|v|_{j,p,T}
\le {d+j \choose j}^{\frac{1}{p}}\tilde{C}(j,p,d,\sigma_0)
\left( \frac{h_T}{k^2} \right)^{-j}
\left(\frac{h_T}{ 2 (q+1)k^2} \right)^{d(\frac{1}{p}-\frac{1}{q})} \|v\|_{0,q,T}.
$$
Under the assumption
$h_T^{-1} \geq 1$ for $0 \le j \le l,$
the corresponding Sobolev norm can be estimated as
\begin{equation}
\begin{array}{rcl}
\|v\|_{j,p,T}^p &=& \sum_{s=0}^j |v|_{s,p,T}^p \\
&\le& \sum_{s=0}^j {d+s \choose s}\left\{ \tilde{C}(s,p,d,\sigma_0)
\left( \frac{h_T}{k^2} \right)^{-s} \left(\frac{h_T}{ 2 (q+1)k^2} \right)^{d(\frac{1}{p}
-\frac{1}{q})} \right\}^p \|v\|_{0,q,T}^p \\
&\le& {d+1+j \choose j} \left\{ \tilde{C}(j,p,d,\sigma_0)
\left( \frac{h_T}{k^2} \right)^{-j} \left(\frac{h_T}{ 2 (q+1)k^2} \right)^{d(\frac{1}{p}
-\frac{1}{q})} \right\}^p \|v\|_{0,q,T}^p,
\end{array}
\label{inv_m_0}
\end{equation}
and the assertion is proved for $m=0.$
Now, let $0 \le m \le l$ and $\alpha$ be a multiindex with $0 \le |\alpha|\le l.$
In the case $|\alpha| \le l-m,$ the estimate $(\ref{inv_m_0})$ and the relation
$$
\|\partial^\alpha v \|_{0,p,T}^p
\le \sum_{|\alpha|\le l-m} \|\partial^\alpha v\|_{0,p,T}^p=\|v\|_{l-m,p,T}^p
$$
imply that
$$
\begin{array}{rcl}
\|\partial^\alpha v\|_{0,p,T}
&\le& C(l-m,p,d,\sigma_0) \left( \frac{h_T}{k^2} \right)^{m-l}
\left(\frac{h_T}{ 2 (q+1)k^2} \right)^{d(\frac{1}{p}-\frac{1}{q})}\|v\|_{0,q,T} \\
&\le& C(l-m,p,d,\sigma_0) \left( \frac{h_T}{k^2} \right)^{m-l}
\left(\frac{h_T}{ 2 (q+1)k^2} \right)^{d(\frac{1}{p}-\frac{1}{q})}\|v\|_{m,q,T},
\end{array}
$$
where $C(j,p,d,\sigma_0):={d+1+j \choose j}^{1/p} \tilde{C}(j,p,d,\sigma_0).$
This inequality is also valid for $l-m \le |\alpha|\le l,$
because there exist two further multiindices $\beta$ and $\gamma$ such that
$\alpha=\beta+\gamma,$ $|\beta|=l-m$ and $|\gamma|\le m.$
As in the first case we get
$$
\begin{array}{rcl}
\|\partial^\alpha v\|_{0,p,T}
&=&\|\partial^\beta \left( \partial^\gamma v \right) \|_{0,p,T} \\
&\le& C(l-m,p,d,\sigma_0) \left( \frac{h_T}{k^2} \right)^{m-l}
\left(\frac{h_T}{ 2 (q+1)k^2} \right)^{d(\frac{1}{p}-\frac{1}{q})}
\|\partial^\gamma v\|_{0,q,T} \\
&\le& C(l-m,p,d,\sigma_0) \left( \frac{h_T}{k^2} \right)^{m-l}
\left(\frac{h_T}{ 2 (q+1)k^2} \right)^{d(\frac{1}{p}-\frac{1}{q})} \|v\|_{m,q,T}.
\end{array}
$$
Since this estimate is valid for $0\le |\alpha|\le l,$
the norm $\|v\|_{l,p,T}^p$ can be estimated by using its definition,
and the statement follows with
$$
C(l,m,p,d,\sigma_0):=C_{l-m,d}\left(2 \sqrt{d} \sigma_0 C_{M,p} \right)^{l-m}
{d+1+l-m \choose l-m}^{1/p} {d+l \choose l}^{1/p}.
$$
The equality sign occurs in the case $l=m=0$ and $p=q.$
\close

\begin{lemma}[global inverse inequality] %
Under the assumptions of Lemma $\ref{inverse_estimate},$
for a locally quasiuniform family of affine partitions
$\{\mathcal{T}_h\}_{h\in(0,1]}$ and for all $v \in W_h,$ there exist
\begin{enumerate}
\item
for $1 \le p,q \le \infty$ a constant
$C=C(l,m,p,q,d,\sigma_0,\hat{T},P(\hat{T}),C_{qu})>0$
such that
$$
\Big( \sum_{T \in \mathcal{T}_h} \|v\|_{l,p,T}^p \Big)^{\frac{1}{p}}
\le C h^{m-l+min(0,d(\frac{1}{p}-\frac{1}{q}))}
\Big( \sum_{T \in \mathcal{T}_h}\|v\|_{m,q,T}^q \Big)^{\frac{1}{q}},
$$
\item
for $1 \le q \le p \le \infty$ and %
$\hat{P}=\mathbb{Q}_k(\hat{T})$ a constant $C=C(l,m,p,d,\sigma_0)>0$
$$
\Big( \sum_{T \in \mathcal{T}_h} \|v\|_{l,p,T}^p \Big)^{\frac{1}{p}}%
\le C \left(\frac{C_{qu} h}{k^2}\right)^{m-l}
\left( \frac{C_{qu} h}{2 (q+1)k^2} \right)^{d(\frac{1}{p}-\frac{1}{q})}
\Big( \sum_{T \in \mathcal{T}_h} \|v\|_{m,q,T}^q \Big)^{\frac{1}{q}}.
$$
\end{enumerate}
\end{lemma}
\proof
The proof of the first estimate can be found in
\cite[Cor.\ 1.141]{Ern:04}.
With
$$
C=C(l,m,p,d,\sigma_0),
$$
the second estimate follows from Lemma \ref{inverse_estimate}
and from the quasiuniformity of $\{\mathcal{T}_h\}_{h>0}$:
$$
\sum_{T \in \mathcal{T}_h} \|v\|_{l,p,T}^p
\le C^p \left(\frac{C_{qu}h}{k^2}\right)^{p(m-l)}
\left( \frac{C_{qu}h}{2 (q+1)k^2} \right)^{pd(\frac{1}{p}-\frac{1}{q})}
\sum_{T \in \mathcal{T}_h} \|v\|_{m,q,T}^p.
$$
Taking the $p$-th root, together with
$\|\cdot\|_{\ell^p} \le \|\cdot\|_{\ell^q}$
the statement follows for $p,q\ne \infty.$
The proof in the case $p=q=\infty$ immediately follows from the
definitions of the norms.
\close

\small
\bibliographystyle{alpha}

\end{document}